\documentclass[12pt]{article}
\usepackage{amsmath,amssymb}
\usepackage{graphicx}
\usepackage[T1, T2A]{fontenc}
\def\date{\hfill \number\day.  \number\month. \number\year }
\let\tilde=\widetilde

\def\ledot{\le\!\!\!\raise2pt\hbox{$\scriptscriptstyle\bullet$}\!}
\def\ldot{\!<\!\!\!\raise1.5pt\hbox{$\scriptscriptstyle\bullet$}}
\def\and{\hskip2pt{\scriptstyle \land}\hskip2pt}

\def\smcap{\kern 2pt {\scriptstyle \cap}\kern 2pt}
\def\smcup{\kern 2pt {\scriptstyle \cup}\kern 2pt}
\def\ssm{{\smallsetminus}}

\def\Qed{\hglue 0pt plus 1filll $\scriptstyle\square$}  

\def\cB{{\cal B}}
\def\cC{{\cal C}}
\def\cD{{\cal D}}
\def\cE{{\cal E}}
\def\cF{{\cal F}}
\def\cG{{\cal G}}
\def\cH{{\cal H}}

\def\cM{{\cal M}}

\def\cO{{\cal O}}
\def\cP{{\cal P}}

\def\1{1\kern-3.pt {\rm l}}    

\def\Cs#1#2{{\rm Cs}_{\hskip1pt #1}{\hskip.5pt#2}} 

\def\Aut{\mathop{{\rm Aut}}}

\def\SL#1#2{{\rm SL}_{#1}{#2}}  
\def\PSL#1#2{{\rm PSL}_{#1}{#2}} 

\def\Opr#1#2{{\rm O}^\prime_{#1}{#2}}  

\def\SO#1#2{{\rm SO}_{#1}{#2}}

\def\Spin#1#2{{\rm Spin}_{#1}{#2}}   
\def\U#1#2{{\rm U}_{#1}{#2}}
\def\SU#1#2{{\rm SU}_{#1}{#2}} 
\def\PU#1#2{{\rm PU}_{#1}{#2}}
\def\PSU#1#2{{\rm PSU}_{#1}{#2}}  
\def\Sp#1#2{{\rm Sp}_{#1}{#2}}
\def\PSp#1#2{{\rm PSp}_{#1}{#2}}

\def\Gtwo{{\rm G}_2}

\font\Bbb=msbm10 at12 pt
\font\scriptBbb=msbm7  at16 pt
\font\scriptscriptBbb=msbm5 at12 pt

\newfam\Bbbfam
\textfont\Bbbfam=\Bbb
\scriptfont\Bbbfam=\scriptBbb
\scriptscriptfont\Bbbfam=\scriptscriptBbb

\def\CC{{\fam=\Bbbfam C}}

\def\HH{{\fam=\Bbbfam H}}

\def\OO{{\fam=\Bbbfam O}}

\def\RR{{\fam=\Bbbfam R}}
\def\sS{{\fam=\Bbbfam S}}

\def\ZZ{{\fam=\Bbbfam Z}}

\font\sansserif=cmss10 at 12 pt 
\font\scriptsansserif=cmss10 at 9 pt
\font\scriptscriptsansserif=cmss10 at 6 pt
\newfam\ssfam
\textfont\ssfam=\sansserif
\scriptfont\ssfam=\scriptsansserif
\scriptscriptfont\ssfam=\scriptscriptsansserif

\let\usuDelta=\Delta  
\let\usuGamma=\Gamma  
\let\usuLambda=\Lambda  
\let\usuOmega=\Omega  
\let\usuPhi=\Phi  
\let\usuPi=\Pi  
\let\usuPsi=\Psi  
\let\usuSigma=\Sigma  
\let\usuTheta=\Theta  
\let\usuUpsilon=\Upsilon  
\let\usuXi=\Xi  

\def\Alpha{{\fam=\ssfam A}}
\def\Beta{{\fam=\ssfam B}}
\def\Chi{{\fam=\ssfam \usuChi}}
\def\Delta{{\fam=\ssfam \usuDelta}}

\def\Eta{{\fam=\ssfam H}}
\def\Gamma{{\fam=\ssfam \usuGamma}}
\def\Chi{{\fam=\ssfam X}}

\def\Kappa{{\fam=\ssfam K}}
\def\Lambda{{\fam=\ssfam \usuLambda}}

\def\Nu{{\fam=\ssfam N}}
\def\Omega{{\fam=\ssfam \usuOmega}}
\def\Phi{{\fam=\ssfam \usuPhi}}
\def\Pi{{\fam=\ssfam \usuPi}}
\def\Psi{{\fam=\ssfam \usuPsi}}

\def\Rho{{\fam=\ssfam P}}
\def\Sigma{{\fam=\ssfam \usuSigma}}
\def\Tau{{\fam=\ssfam T}}
\def\Theta{{\fam=\ssfam \usuTheta}}
\def\Upsilon{{\fam=\ssfam \usuUpsilon}}
\def\Ypsilon{{\fam=\ssfam \usuUpsilon}}
\def\Xi{{\fam=\ssfam \usuXi}}
\def\Zeta{{\fam=\ssfam Z}}

\let\epsilon=\varepsilon
\let\theta=\vartheta
\let\phi=\varphi
\let\rho=\varrho

\chardef\tempcat=\the\catcode`\@
\catcode`\@=11

\def\cydot{{\mathsurround=0pt$\cdot$}}

\def\ubar#1{\oalign{#1\crcr\hidewidth
    \vbox to.2ex{\hbox{\char22}\vss}\hidewidth}}

\def\cprime{\/{\mathsurround=0pt$'$}}
\def\Cprime{{\mathsurround=0pt$'$}}
\def\cdprime{\/{\mathsurround=0pt$''$}}
\def\Cdprime{{\mathsurround=0pt$\ubar{\hbox{$''$}}$}}

\def\dbar{dj}           
\def\Dbar{Dj}           

\def\dz{dz}
\def\Dz{Dz}
\def\dzh{dzh\cydot }
\def\Dzh{Dzh\cydot }

\def\@gobble#1{}
\def\@testgrave{\`}
\def\@stressit{\futurelet\chartest\@stresschar }

\def\@stresschar#1{%
  \ifx #1y\def\result{\futurelet\chartest\@yligature}%
  \else \ifx #1Y\def\result{\futurelet\chartest\@Yligature}%
  \else \ifx\chartest\@testgrave \def\result{\accent"26 }%
  \else \def\result{\accent"26 #1}%
  \fi \fi \fi
  \result }

\def\@yligature{%
  \ifx a\chartest \def\result{\accent"26 \char"1F \@gobble}%
  \else \ifx u\chartest \def\result{\accent"26 \char"18 \@gobble}%
  \else \def\result{\accent"26 y}%
  \fi \fi
  \result }

\def\@Yligature{%
  \ifx a\chartest \def\result{\accent"26 \char"17 \@gobble}%
  \else \ifx A\chartest \def\result{\accent"26 \char"17 \@gobble}%
  \else \ifx u\chartest \def\result{\accent"26 \char"10 \@gobble}%
  \else \ifx U\chartest \def\result{\accent"26 \char"10 \@gobble}%
  \else \def\result{\accent"26 Y}%
  \fi \fi \fi \fi
  \result }

\def\!{\ifmmode \mskip-\thinmuskip \fi}

\def\cyracc{%
  \def\cydot{{\kern0pt}}%
  \def\cprime{\char"7E }\def\Cprime{\char"5E }%
  \def\cdprime{\char"7F }\def\Cdprime{\char"5F }%
  \def\dbar{dj}\def\Dbar{Dj}%
  \def\dz{\char"1E }\def\Dz{\char"16 }%
  \def\dzh{\char"0A }\def\Dzh{\char"02 }%
  \def\'##1{\if c##1\char"0F %
    \else \if C##1\char"07 %
    \else \accent"26 ##1\fi \fi }%
  \def\`##1{\if e##1\char"0B %
    \else \if E##1\char"03 %
    \else \errmessage{accent \string\` not defined in cyrillic}%
        ##1\fi \fi }%
  \def\=##1{\if e##1\char"0D %
    \else \if E##1\char"05 %
    \else \if \i##1\char"0C %
    \else \if I##1\char"04 %
    \else \errmessage{accent \string\= not defined in cyrillic}%
        ##1\fi \fi \fi \fi }%
  \def\u##1{\if \i##1\accent"24 i%
    \else \accent"24 ##1\fi }%
  \def\"##1{\if \i##1\accent"20 \char"3D %
    \else \if I##1\accent"20 \char"04 %
    \else \accent"20 ##1\fi \fi }%
  \def\!{\ifmmode \def\result{\mskip-\thinmuskip}%
    \else \def\result{\@stressit}\fi \result}}

\def\keep@cyracc{\let\cyr=\relax \let\i=\relax
        \let\ubar=\relax \let\cydot=\relax
        \let\cprime=\relax \let\Cprime=\relax
        \let\cdprime=\relax \let\Cdprime=\relax
        \let\dbar=\relax \let\Dbar=\relax
        \let\dz=\relax \let\Dz=\relax
        \let\dzh=\relax \let\Dzh=\relax
        \let\'=\relax \let\`=\relax \let\==\relax
        \let\u=\relax \let\"=\relax \let\!=\relax }

\catcode`\@=\tempcat

\font\tencyr=wncyr10 scaled1200

\def\cyr{\tencyr\cyracc}
\font\tencyss=wncyss10 scaled1200
\def\cyss{\tencyss\cyracc}

\def\sm{\hskip-1pt$\raise1pt\hbox{$\scriptstyle\setminus$}$\hskip-1pt}

\begin{document}
\abovedisplayskip=3pt
\belowdisplayskip=3pt

\overfullrule=.1pt
\font\bf=cmbx10 scaled 1200
\font\Bf=cmbx10 scaled 1500
\font\bbf=cmbx10 scaled 1800
\font\sbf=cmbx12 scaled 970
\font\ss=cmss12
\def\bullett{\raise1pt\hbox{$\scriptscriptstyle\bullet$}} 
\def\ssm{\smallsetminus}
\def\rk{\rm rk\,}
\let\hat=\widehat

\let\bold=\bf
\def\bib{\bibitem{} }
\def\Qed{\hglue 0pt plus 1filll $\square$}
\font\teneufm=eufm10 scaled 1300 
\font\seveneufm=eufm10 scaled 1000
\newfam\eufmfam
\textfont\eufmfam=\teneufm
\scriptfont\eufmfam=\seveneufm

\def\frak#1{{\fam\eufmfam\relax#1}}
\def\frL{{\frak L}}
\def\fre{{\frak e}}

\def\3ast{$(\lower1pt\hbox{$\scriptstyle**$})\hskip-11pt\raise3.2pt\hbox{$\scriptstyle*$}\hskip6pt$}
\def\ldot{\!<\!\!\!\raise1.3pt\hbox{$\scriptscriptstyle\bullet$}}
\def\Rtimes{\times\hskip-3pt\raise0pt\hbox{$\vrule height 6pt width .8pt$}}
\def\smotimes{\begin{scriptsize}$\otimes$\end{scriptsize}}
\def\smoplus{\begin{scriptsize}$\oplus$\end{scriptsize}} 
\def\tsmotimes{\begin{scriptsize}$\tilde\otimes$\end{scriptsize}}
\def\tsmoplus{\begin{scriptsize}$\tilde\oplus$\end{scriptsize}} 
\def\GGamma{\Gamma\hskip-4.6pt\Gamma}
\def\Yu{{\cyss Yu}}

\centerline{\bbf Semi-simple groups of  } 
\par\smallskip
\centerline{\bbf compact 16-dimensional planes}
\par\bigskip
\centerline {by Helmut R. \sc Salzmann}  
\par\bigskip
\begin{abstract}
The automorphism group $\Sigma$ of a compact topological projective plane with a $16$-dimensional point space is a locally compact group. If the dimension of $\Sigma$ is at least $29$, then $\Sigma$ is known to be a Lie group. For the connected component $\Delta$ of $\Sigma$ the condition 
$\dim\Delta{\,\ge\,}27$ suffices. Depending on the structure of $\Delta$ and the configuration of the fixed elements of $\Delta$ sharper bounds are obtained here. Example: If $\Delta$ is semi-simple and fixes exactly one line and possibly several points on this line, then $\Delta$ is a Lie group if $\dim\Delta{\,\ge\,}11$. Any semi-simple group 
$\Delta$ which satisfies $\dim\Delta{\,\ge\,}25$ is a 
Lie group.
\end{abstract}
\par\bigskip

{\Bf 1. Introduction}
\par\medskip
A compact connected topological projective plane has a point space $P$ 
of (covering) dimension $\dim P{\,=\,}d{\,|\,}16$, provided $d{\,<\,}\infty$. 
For properties and resuts concerning $16$-dimensional planes see 
the treatise \cite{cp} and a more recent update \cite{66}. 
With the compact-open topology (the topology of uniform convergence), the 
automorphism group $\Sigma{\,=\,}\Aut\cP$ of a compact $16$-dimensional plane $\cP{\,=\,}(P,\frL)$ is a locally compact transformation group of $P$ 
as well as of the line space $\frL$, see \cite{cp} 44.3. 
In the case of the {\it classical\/} plane $\cO$ over the octonion algebra 
$\OO$ the group $\Sigma$ is a simple Lie group of 
dimension $78$\; (\cite{cp} 18.19). If $\cP$ is not classical, then
$\dim\Aut\cP{\,\le\,}40$\; (\cite{cp} 87.7). All planes with 
$\dim\Sigma{\,\ge\,}35$ are known explicitly, provided $\Sigma$ does 
not fix exactly one incident point-line pair (\cite{66}). 
{\it If $\dim\Sigma{\,\ge\,}27$, then the connected component $\Delta$ 
of $\Sigma$ is a Lie group\/}\;  (\cite{psz}). More detailed results can be found in \cite{66}. For semi-simple groups 
sharper bounds will be obtained here.

\par\bigskip
{\Bf 2. Preliminaries} 
\par\medskip
In the following, $\cP{\,=\,}(P,\frL)$ will always mean a compact 
$16$-dimensionl projective plane if not stated otherwise; $\Delta$ denotes 
a connected closed subgroup of $\Aut\cP$.
\par\medskip
{\bf Notation} is more or less standard and agrees with that in the book \cite{cp}. A \emph{flag} is an incident point-line pair, a \emph{double flag} consists of two points, say $u,v$, their join $uv$, and a second line in the \emph{pencil}~$\frL_v$. A $2$-dimensional plane will also be called {\it flat\/}. Homeomorphism is indicated by $\approx$. The topological dimension 
{\it dim\/} of a set $M{\subseteq}P$ or a group, the covering dimension, coincides with the inductive dimension, see \cite{cp} \S\hskip1pt92,\,93.5.
As  customary, $\Cs\Delta\Gamma$ or just $\Cs{}\Gamma$ is the centralizer of $\Gamma$ in $\Delta$. Distinguish between the commutator subgroup 
$\Gamma'$ and the connected component $\Gamma^1$ of the 
topological  group $\Gamma$. Local isomorphy of groups is symbolized by 
$\circeq$. The coset space 
$\Delta/\Gamma{\,=\,}\{\Gamma\delta\mid\delta{\,\in\,}\Delta\}$ has the dimension 
$\Delta{:}\Gamma{\hskip1pt=\hskip1pt}\dim\Delta{-}\dim\Gamma$.
If $\Gamma$ is any subset of $\Delta$, then $\cF_\Gamma$ denotes 
the configuration of all fixed elements (points and lines) of $\Delta$.
The group $\Delta_{[c,A]}$ consists of the axial collineations  in $\Delta$ 
with axis $A$ and  center~$c$. A collineation group $\Gamma$ is said to
be \emph{straight} if each orbit $x^\Gamma$ is contained in some line. 
In this case a theorem of Baer \cite{ba} asserts that either 
$\Gamma{=\,}\Gamma_{\hskip-2pt[c,A]}$ is a group of axial collineations  
or  $\cF_\Gamma$ is a Baer subplane\, (cf. 2.2 below).
\par\medskip
{\bf 2.1 Topology.}   {\it Each line $L$ of $\cP$  is homotopy equivalent $(\simeq)$ to an $8$-sphere $\sS_8$\/}, see  \cite{cp} 54.11. So far, no example with $L{\,\not\approx\,}\sS_8$ has been found. 
The Lefschetz fixed point theorem implies that  each homeomorphism $\phi{\,:\,}P{\,\to\,}P$ has a fixed point. By duality, {\it each automorphism of $\cP$ fixes a point and a line\/}, 
see \cite{cp} 55.\,19,\,45. 
\par\medskip
{\bf 2.2 Baer subplanes.} {\it Each $8$-dimensional $closed$ subplane 
$\cB$ of   $\cP$ is a Baer subplane\/}, i.e., each point of $\cP$ is incident 
with a line of $\cB$ (and dually, each line of $\cP$ contains a point of 
$\cB$), see \cite{57} \S~\hskip-3pt3 or \cite{cp} 55.5 for details.  
By a theorem of L\"owen \cite{lw} Th.\,B, any two Baer subplanes of a compact 
connected projective plane of finite dimension have a point and a line in common.
$\langle \cM\rangle$ will denote the smallest 
{\it closed\/} subplane of $\cP$ containing the  set $\cM$ of points and lines. We write $\cB{\,\ldot}\cP$ if $\cB$ is a Baer subplane and $\cB{\,\ledot\,}\cP$ if $\,8\vert\dim\cB$. 
\par\medskip
{\bf 2.3 Groups.} {\it Any {\rm connected} subgroup $\Delta$ of 
$\Sigma{\,=\,}\Aut\cP$ with $\dim\Delta{\,\ge\,}27$ is a Lie  group\/}, see \cite{psz}. For $\dim\Delta{\,\ge\,}35$, the result has been proved in \cite{cp} 87.1.  In particular, $\Delta$ is then either semi-simple, or $\Delta$ has a central torus subgroup or a minimal normal vector  subgroup 
$\Xi{\,\cong\,}\RR^k$, see \cite{cp} 94.26.
  
\par\medskip
{\bf 2.4 Involutions.} Each involution $\iota$ is either a reflection or it is {\it planar\/} 
($\cF_{\hskip-1.5pt\iota}{\,\ldot}\cP$), see \cite{cp} 55.29. Commuting involutions with the same fixed point set are identical (\cite{cp} 55.32). Let 
$\ZZ_2^{\hskip3pt r}{\,\cong\,}\Phi{\,\le}\Aut\cP$ and 
$\dim P{\,=\,}2^m$.  Then $r{\,\le\,}m{+}1$; if $\Phi$ is generated by reflections, then $r{\,\le\,}2$, see \cite{cp} 55.34(c,b). 
If $r{\,\ge\,}m \ ({\,>\,}2)$, then  $\Phi$ contains a  reflection (a planar 
involution), cf. \cite{cp} 55.34(d,b). Any torus group in $\Aut\cP$ has dimension at most $m$,\, see \cite{cp} 55.37, for $m{\,=\,}4$ in particular, $\rk\Delta{\,\le\,}4$,  see also~\cite{cs}. The orthogonal group $\SO5\RR$ cannot act non-trivially on $\cP\,$ (\cite{cp} 55.40). 
\par\medskip
{\bf 2.5 Stiffness} refers to the fact that the dimension of the stabilizer 
$\Lambda$ of a (non-degenerate) quadrangle $\fre$ cannot be very large.  
If $\cP$ is the classical plane, then   $\Lambda$ is  the $14$-dimensional exceptional compact simple Lie group $\Gtwo{\,\cong\,}\Aut\OO$ for 
{\it each\/} choice of $\fre$, cf. \cite{cp} 11.\,30--35. In the general case, let $\Lambda$ denote the connected component of  $\Delta_\fre$. Then  the following holds: \par
\quad (a) {\it $\Lambda{\,\cong\,}\Gtwo$ and $\cF_\Lambda$ is flat, or $\dim\Lambda{\,<\,}14$\/} 
  (\cite{cp} 83.\,23,\,24), \par
\quad (\^a) {\it $\Lambda{\,\cong\,}\Gtwo$ or $\dim\Lambda{\,\le\,}11$\/}\ (\cite{B1} 4.1), \par
\quad (b) {\it if $\cF_\Lambda$ is a Baer subplane, then $\Delta_\fre$ is compact and  
 $\dim\Lambda{\,\le\,}7$\/} \  (\cite{cp} 83.6), \par
\quad (\^b) {\it if $\cF_\Lambda{\;\ldot}\cP$ and $\Lambda$ is a Lie group, then  $\Lambda{\,\cong\,}\Spin3\RR$ or  
 $\dim\Lambda{\,\le\,}1$\/} \ (\cite{cp} 83.22), \par
\quad (c) {\it if there exists a subplane $\cB$ such that 
   $\cF_\Lambda{\;\ldot}\cB{\;\ldot}\cP$, then $\Lambda$ is compact\/} \par\hskip30pt   (\cite{37} 2.2, \cite{cp} 83.9), \par
\quad (\^c) {\it if $\Lambda$ contains a pair of commuting invo lutions, then  $\Delta_\fre$ is 
 compact\/} \ (\cite{cp} 83.10), \par
\quad (d) {\it if $\Lambda$ is compact or semi-simple or if  $\cF_\Lambda$   is connected, then 
 $\Lambda{\,\cong\,}\Gtwo$ or $\dim\Lambda{\,\le\,}10$\/} \   \par\hskip30pt (\cite{gs} XI.9.8, \cite{37} 4.1, \cite{B1} 3.5), \par
\quad (e) {\it if $\Lambda$ is compact, or if 
$\dim\cF_\Lambda{\,=\,}4$, or if    $\Lambda$ is a Lie group and $\cF_\Lambda$ is connected,   then \par\hskip30pt
$\Lambda{\,\cong\,}\Gtwo$ or $\Lambda{\,\cong\,}\SU3\CC$ or $\dim\Lambda{\,<\,}8$\/} \
   (\cite{37} 2.1; \cite{44}, \cite{B1} 3.5; \cite{B2}), \par
\quad (\^e) {\it if $\Lambda$ is a compact Lie group, then 
$\Lambda{\hskip1pt\cong\,}\Gtwo$,  
$\SU3\CC,\; \SO4\RR$,  or $\dim\Lambda{\,\le\,}4$\/}\,\ (\cite{37} 2.1). 
\par\medskip

{\bold 2.6 Groups with open orbits.}
{\it Let $X$ denote the point space $P$ of $\cP$ or a line
$L$ of $\cP$. If $U{\,\subseteq\,}X$ is a $\Delta$-orbit which is open in $X$ or, equivalently, satisfies $\dim U{\,=\,}\dim X$, then $X$ is a manifold and $\Delta$ induces a Lie group on $U$. If $X{\,=\,}L$, then it follows that all lines are mani\-folds homeomorphic to $\sS_8$. If $X{\,=\,}P$, then $\Sigma$ is a Lie group\/}.\, (\cite{cp} 53.2).
\par\smallskip
{\tt Remark.} The result  \cite{cp} 53.2 depends on a theorem of Szenthe 
for which a correct proof has been given in the meantime by Hofmann 
and Kramer \cite{HK} 5.5 Corollary.

\par\medskip
{\bf 2.7 Lemma.} {\it  If $\cF_{\hskip-2pt\Delta}{\,=\,}\emptyset$, then $\Delta$ is semi-simple with trivial center, or $\Delta$ induces a centerfree semi-simple group on some connected proper closed subplane. In particular, 
$\cF_{\hskip-2pt\Delta}{\,\ne\,}\emptyset$ whenever $\Delta$ is 
commutative\/}. 
\par\smallskip
{\tt Proof.} As a flat plane does not have a proper closed subplane, 
the Lemma is true for flat planes by \cite{cp} 33.1, and each commutative 
connected subgroup fixes a point or a line. Suppose the Lemma has been shown for all planes of smaller dimension thane the plane on which 
$\Delta$ acts. If $\Delta$ contains a central element $\zeta{\,\ne\,}\1$ or
 if $\Delta$ is not semi-simple and hence has a commutative normal subgroup $\Xi$, then there is a point $p$ such that $p^\zeta{\,=\,}p$ or
$p^\Xi{\,=\,}p$ (up to duality). By assuption, $p^\Delta{\,\ne\,}p$ and
$p^\Delta$ is not contained in a line. Obviously $\zeta|_{p^\Delta}{\,=\,}\1$,
and normality of $\Xi$ implies $\Xi|_{p^\Delta}{\,=\,}\1$. Therefore 
$\cE{\,=\,}\langle p^\Delta\rangle$ is a proper connected subplane 
and $\overline\Delta{\,=\,}\Delta|_\cE$ satisfies 
$\cF_{\hskip-1pt\overline\Delta}{\,=\,}\emptyset$. Now the Lemma can be applied to $\overline\Delta$. \Qed
\par\smallskip
{\tt Note} that the Lemma holds for all compact connected planes of finite dimension.
\par\medskip
{\bf 2.8 Addendum.} {\it If $\cF_\Delta{\,=\,}\emptyset$ and if there is no 
$\Delta$-invariant proper subplane, then each involution in a proper simple factor of 
$\Delta$ is  \emph{planar}. This is true for each compact plane of dimension 
$d\,|\,16$\/}. 
\par\smallskip
{\tt Proof.} By 2.7 the group $\Delta$ is semi-simple with trivial center. 
Suppose that $\Delta{\,=\,}\Gamma{\times}\Psi$ is a product of two 
proper factors.  If $\Gamma$ contains a reflection 
$\gamma{\,\in\,}\Gamma_{\hskip-2pt[c,A]}$, then 
$c^\Psi{=\,}c{\,\ne\,}c^\Gamma$ and $\Psi|_{c^\Gamma}{\,=\,}\1$. 
A line $L$ containing $c^\Gamma$ would be fixed by 
$\Gamma\Psi{\,=\,}\Delta$. Therefore 
$\langle c^\Gamma\rangle{\,=\,}\cE{\,\le\,}\cF_\Psi$, and 
$\cE^\Gamma{=\,}c^\Psi{ =\,}c^\Delta$ contrary to the assumtion. 
Hence each involution in a proper factor is {planar}. 
\par\smallskip
{\tt Remark.} The statement is not true if there is an invariant subplane: 
in the  classical quaternion plane $\cH$, the group $\Delta$ consisting of all 
collineations which map the real subplane onto itself is the direct product 
of the automorphism group $\Gamma{\,\cong\,}\SO3\RR$ of the quaternions 
and the group $\Rho{\,\cong\,}\SL3\RR$ of real collineations. The involutions in $\Rho$ are reflections of $\cH$.
\par\medskip 
{\bf 2.9 Fact.} {\it If $\cE^\Delta{\,=\,}\cE{\,<\,}\cP$, if 
$\Delta|_\cE{\,=\,}\Delta/\Kappa$ is a semi-simple Lie group, and if 
$\Lambda{\,=\,}\Kappa^1$, then $\Delta$ contains a covering group 
$\Gamma$ of $\Delta/\Kappa$ such that 
$\Gamma|_\cE{\,\cong\,}\Delta/\Kappa$ and 
$\Delta{\,=\,}\Gamma\Lambda$. In particular, $\Delta$ is a Lie group 
whenever $\Lambda$ is\/}.
\par\smallskip
{\tt Proof.} The first claim is a special case of \cite{cp} 94.27; the second follows from 
$\dim\Gamma{\,=\,}\Delta{:}\Kappa=\dim\Delta{\,-\,}\dim\Lambda$ 
because $\Delta$ is connected. \Qed
\par\bigskip
{\Bf 3. No fixed elements}
\par\medskip
The following slightly improves results in \cite{66}: 
\par\smallskip
{\bold 3.1.} {\it If $\Delta$ is semi-simple, 
$\cF_{\!\Delta}{\,=\,}\emptyset$ and if  $\dim\Delta{\,\ge\,}20$, then $\Delta$ is a Lie group.  For $\dim\Delta{\,\ge\,}21$ the assertion has been proved by different methods in \cite{66} 3.0.}
\par\smallskip
{\tt Proof.} Suppose that $\Delta$ is not a Lie group. \\
(a) By the Lemma above, $\Delta$ is a direct product of simple groups 
with trivial center, or it induces such a group on some $\Delta$-invariant 
connected  subplane $\cE$. Choose $\cE$ in such a way that the dimension 
of $\cE$ is minimal.
In the case $\cE{\,=\,}\cP$ the center of $\Delta$
is trivial and $\Delta$ is a Lie group. If $\cE$ is flat, then
$\Delta|_\cE{\,=\,}\Delta/\Kappa$ is strictly simple of dimension  $3$ or
$\Delta|_\cE{\,\cong\,}\SL3\RR$\;  (cf. \cite{cp} 33.6,7); by Stiffness,
$\Lambda{\,=\,}\Kappa^1{\,\cong\,}\Gtwo$ is a Lie group and 2.9 applies,  or 
$\dim\Kappa{\,\le\,}10$ and $\dim\Delta{\,\le\,}18$. \\
(b) Next, assume that $\cE{\,\ldot}\cP$. All almost simple groups of dimension 
at least $10$ acting on $\cE$ have been determined by Stroppel~\cite{str}; for those satisfying $\cF_\Delta{\,=\,}\emptyset$ see also \cite{65} 2.1. 
Suppose first that 
$\Delta^{\!*}{=\,}\Delta|_\cE{\,=\,}\Delta/\Kappa{\,=\,}\Gamma{\times}\Psi$ 
is not simple, and that $\Gamma$ is a simple factor. By the addendum above, each involution in $\Gamma$ is planar. As the center of 
$\Gamma$ is trivial, there are distinct involutions 
$\alpha,\beta{\,\in\,}\Gamma$. 
Because of L\"owen's theorem (see 2.2 above)  there is some point 
$c{\,\in}\cF_{\alpha,\beta}$, and the point space $F_{\alpha,\beta}$ of 
$\cF_{\alpha,\beta}$  has dimension at most~$2$. If $c^\Psi{=\,}c$, then 
$\langle c^\Delta\rangle{\,=\,}\langle c^\Gamma\rangle{\,\le\,}\cF_\Psi{\,<\,}\cE$ contrary to the assumption. Similarly, if $c^\Psi{\ne\,}c$ is contained in a line $L$ of $\cE$, then $L{\,=\,}L^\Psi$ induces a line of $\cF_\alpha$ and of $\cF_\beta$, 
a contradiction dual to the previous one. Hence $\langle c^\Psi\rangle$ is a 
flat subplane, and $\Psi$ is simple of dimension $3$ or $8$. 
Moreover, $\Lambda{\,=\,}\Kappa^1$ is a compact normal, hence semi-simple subgroup, in fact a Lie group (see \cite{cp} 93.11). Stiffness implies 
$\dim\Lambda{\,\le\,}6$.     
Interchanging the r\^oles of $\Gamma$ and $\Psi$, it follows that both
$\Gamma$ and $\Psi$ have dimension $8$, in fact, 
$\Gamma{\,\cong\,}\Psi{\,\cong\,}\SL3\RR$, and $\Delta$ is a Lie group by 2.9. \\ 
(c) According to \cite{65} 2.1,2, the induced group $\Delta^{\!*}{=\,}\Delta|_\cE$ is a motion group of the classical quaternion plane or $\Delta^{\!*}{\,\cong\,}\SL3\CC$ acts on a Hughes plane and  fixes a proper subplane of $\cE$. We will have to deal  with this group in the case $\dim\cE{\,=\,}4$. 
If $\Delta^{\!*}{\,\cong\,}\PU3(\HH,r)$, then 
$\dim\Delta^{\!*}{=\,}21$. Recall that $\Delta$ is semi-simple. The kernel 
$\Kappa$ is compact and $\Lambda{\,=\,}\Kappa^1$ is semi-simple; 
$\Delta/\Lambda$  is locally isomorphic to $\U3(\HH,r)$. By \cite{psz} only 
the cases $\dim\Delta{\,<\,}27$ must be considered. Then 
$\dim\Lambda{\,\mid\,}6$ and $\Lambda$ is a compact Lie group. 
For $\Delta^{\!*}{\,\cong\,}\PSU4(\CC,1)$ as a motion group of the
quaternion plane cf. \cite{cp} 18.32. A maximal compact subgroup of 
$\Delta^{\!*}$ is isomorphic to $\U3\CC$ and may be the image of a non-Lie group. 
Again $\Kappa$ is compact by stiffness, it is a finite extension of the semi-simple group $\Lambda{\,=\,}\Kappa^1$, and $\Lambda$ is a Lie group, 
see \cite{cp} 93.11. Thus $\Delta$ is a Lie group. \\
{\tt Note} that in steps (a--c) the smaller bound $\dim\Delta{\,\ge\,}19$ suffices;
the assumption $\dim\Delta{\,\ge\,}20$ will only be needed in the next step. \\
(d) Finally, let $\dim\cE{\,=\,}4$. Suppose that $\dim\Delta|_\cE{\,>\,}3$.
According to \cite{cp} 71.8 the induced group 
$\Delta|_\cE$ is an almost simple Lie group, and   
$\Delta{:}\Kappa{\,\le\,}8$ or $\Delta{:}\Kappa{\,=\,}16$. In this step, the first
 case will be considered. Then $\dim\Kappa{\,\ge\,}12$ and 
$\Lambda{\,=\,}\Kappa^1{\,\cong\,}\Gtwo$. Therefore 
$\Delta{:}\Kappa{\,\ge\,}6$, and 2.9 shows that $\Delta$ is a Lie group. \\ 
(e) In the second case, $\Delta|_\cE{\,\cong\,}\PSL3\CC$ by \cite{cp} 71.8.
This includes the  possibility that there is a Hughes plane 
$\cH^\Delta{\,=\,}\cH$ such that $\cE{\,<\,}\cH{\,\ldot\,}\cP$. 
The plane $\cE$  is the classical complex plane\, (\cite{cp} 72.8), 
and $4{\,\le}\dim\Kappa{\,\le\,}8$. 
The normal  subgroup $\Lambda{\,=\,}\Kappa^1$ is semi-simple and hence 
$\dim\Lambda{\,\in\,}\{6,8\}$.  Consequently 
$\dim\Delta{\,\ge\,}22$.  Choose a line $L$ of $\cE$, a point 
$z{\,\in\,}L\sm\cE$, and two point $a,b$ of $\cE$,\, $a,b{\,\notin\,}L$. 
If $\langle \cE,z\rangle{\,<\,}\cP$, then $\Lambda$ would be a compact Lie 
group by 2.5(c). Hence $\Delta_z$ acts effectively on $\cE$. 
The stabilizer $\Gamma{\,=\,}\Delta_{z,a,b}$ fixes each point of $\cE$ 
on the line $ab$, and $\dim\Gamma{\,\ge\,}3$ because of 2.6.
Let $\Nu$ be a compact central subgroup of $\Lambda$, so that 
$z^\Nu{\,\subset\,}\cF{\,=\,}\cF_\Gamma$. As $\Nu$ is not a Lie group and 
$\Nu_z{\,=\,}\1$, it follows from \cite{cp} 71.2 that $\cF{\,\ldot}\cP$. Now 
$\Gamma$ is compact by 2.5(b). On the other hand, $\Gamma$ is a group 
of axial collineations  in  $\Delta_z$. These form a group
$\CC^2{\times}\CC^{\times}$ and do not contain a compact subgroup of 
dimension~$3$, a contradiction. \Qed
\par\medskip
{\bf 3.2 Theorem.} {\it If $\cF_\Delta{\,=\,}\emptyset$ and 
$\dim\Delta{\,\ge\,}22$, then $\Delta$ is a Lie group\/}. 
\par\smallskip
{\tt Proof.} (a) Suppose that $\Delta$ is not a Lie group. Then 
$\dim\Delta{\,<\,}27$  by \cite{psz}, and there are arbitrarily small compact, $0$-dimensional  central subgroups $\Nu$ such that $\Delta/\Nu$ is  a Lie group, cf. \cite{cp} 93.8. 
Choose a fixed point~$x$  of some element $\zeta{\,\in\,}\Nu\sm\{\1\}$. 
From $\zeta|_{x^\Delta}{\,=\,}\1$ and  $\cF_\Delta{\,=\,}\emptyset$ it follows that $\cD{\,=\,}\langle p^\Delta\rangle$ is a proper connected subplane. 
Put $\Delta|_\cD{\,=\,}\Delta/\Kappa$. If $\cD$ is flat, then 
$\Delta{:}\Kappa{\,=\,}8$ and $\Lambda{\,=\,}\Kappa^1{\,\cong\,}\Gtwo$ 
by stiffness, so that $\Delta|_\cD{\,\cong\,}\SL3\RR$ and 2.9 applies. \\
If  $\dim\cD{\,=\,}4$, then $\Delta/\Kappa$ is a Lie group by 
\cite{cp} 71.2, and stiffness 2.5(e) shows that  $\dim\Kappa{\,<\,}8$
or $\Lambda{\,=\,}\Kappa^1{\,\cong\,}\SU3\CC$.   
Hence $\Delta{:}\Kappa{\,\ge\,}14$, and then 
$\Delta/\Kappa{\,\cong\,}\PSL3\CC$ by \cite{22} 5.6.  
Now $6{\,\le\,}\dim\Lambda{\,\le\,}8$, and the same arguments as at 
the end of step (e) of the previous proof show that 
this is impossible. Consequently $\dim\cD{\,\,=}8$. \\
(b) Let $\cD{\;\ldot}\cP$ and note that $\Delta$ acts on $\cD$ without fixed elements. By  {\cite{65} 2.1} the following holds for 
$\Delta|_\cD{\,=\,}\Delta/\Kappa$:  either  $\cD$ is the classical quaternion plane, $\Delta/\Kappa$ is a motion group of $\cD$ and 
$\Delta{:}\Kappa{\,\in\,}\{21,15\}$,  or $\cD$ is a Hughes plane and
$\Delta'|_\cD{\,\cong\,}\SL3\CC$. Moreover $\Delta/\Kappa$ is a Lie group 
and we may assume that $\Nu|_\cD{\,=\,}\1$. The kernel
$\Kappa$ is compact of dimension at most $7$.  \\
(c)  In the first case, the motion group $\Delta/\Kappa$ is covered by a subgroup 
$\Upsilon{\hskip1pt\circeq\,}\U3(\HH,r)$ of $\Delta$\, (see 2.9), and 
$\dim\Kappa{\,\le\,}5$.  The representation of $\Ypsilon$ on $\Kappa/\Nu$ 
shows that $\Kappa{\,\le\,}\Cs{}\Ypsilon$. We may assume that 
$\Kappa$ is connected, since $\dim\Kappa^1\Ypsilon{\,=\,}\dim\Delta$ and $\Delta$ is connected. Moreover, the center $\Zeta$ of $\Kappa$ has positive dimension: because  $\Kappa$ is not a Lie group, the claim follows from the structure of compact groups as described in \cite{cp} 93.11. The stabilizer $\Pi{\,=\,}\Delta_p$ of a non-absolute point $p{\,\in\,}\cD$ fixes also the polar $L$ of $p$, and $\Pi/\Kappa{\,\circeq\,}\U2\HH{\times}\U1\HH$. In particular, $\Pi$ is compact and $\dim\Pi{\,\ge\,}14$. \\ 
Choose $z{\,\in\,}L\sm\cD$ and note that 
$\Pi_z\smcap\Kappa{\,=\,}\1$ because $\langle \cD,z\rangle{\,=\,}\cP$. Hence $\Pi_z{\,=\,}\Delta_z$  is a Lie group.  From \cite{cp} 96.11 it follows that 
$\Pi{\hskip1pt:\hskip1pt}\Pi_z{\,\le\,}7$ and $\dim\Delta_z{\,\ge\,}7$.  
The group $\Delta_z$ acts faithfully on $\cD$, it is isomorphic to a subgroup 
of $\U2\HH{\times}\U1\HH$. We apply Richardson's theorem \cite{cp} 96.34 
to the action of $\Delta_z$ on $S{\,:=\,}L\smcap\cD{\,\approx\,}\sS_4$\, 
(cf.\cite{cp} 96.22). Put $\Phi{\,=\,}\Delta_z|_S{\,=\,}\Delta_z/\Eta$ and note 
that $\Eta$ is isomorphic to a subgroup of $\U1\HH$, so that $\dim\Phi{\,\ge\,}4$. \\
(d) There are the following possibilities: \\
1) $\Phi{\,\cong\,}\SO3\RR{\times}\SO2\RR$ and $\Phi$ has a one-dimension\-al 
orbit,\; 2) $\Phi$ has two fixed points $u,v{\,\in\,}S$, \\
3) $\Phi{\,\cong\,}\SO5\RR$ is transitive on $S$ and $\dim\Phi{\,\ge\,}10$. 
In the first two cases, there are points $u,v,w{\,\in\,}S$ and $q{\,\in\,}pv$ such 
that $\Delta_{z,u,v,w,q}$ has dimension ${>\,}0$, but this group fixes $S$ 
pointwise and is {trivial} on $\cD$. In the third case, $\Delta_{z,v}$ fixes a second point $u{\,\in\,}S$. Let $p,v{\,\ne\,}q{\,\in\,}pv\smcap\cD{\,:=\,}J$. Then 
$\Gamma{=\,}\Delta_{z,q}$ fixes each point of $J{\,\approx\,}\sS_4$, and 
$\dim\Gamma{\,\ge\,}3$. Hence $\cF_\Gamma{\,\ldot}\cP$, and $J$ is contained 
in a line of $\cF_\Gamma$. By the open mapping theorem \cite{cp} 51.19,20 
we conclude that $J$ is in fact a line of $\cF_\Gamma$. Therefore 
$M{\,=\,}L\smcap\cF_\Gamma{\,\approx\,}\sS_4$, and  
Richardson's theorem  \cite{cp} 96.34 shows that $\Gamma\Kappa$
induces a Lie group on $M$. As $\Kappa_z{\,=\,}\1$, both $\Kappa$ and
$\Delta$ would be Lie groups. \\ 
(e) If $\Delta/\Kappa{\,\cong\,}\PSU4(\CC,1)$ is the $15$-dimensional 
motion group of the planar polarity, then $\Delta$ is doubly transitive on the 
absolute $5$-sphere $U$ in $\cD$ (cf. \cite{cp} 18.32), and $\dim\Kappa{\,=\,}7$. 
In homogeneous coordinates over the quaternion field $\HH$, the planar polarity can be described by the form $\overline xix{+}\overline yiy{+}\overline ziz$, see \cite{cp} 13.18.  The line ``at infinity" $W$, given by  $z{\,=\,}0$, intersects $U$ in the circle consisting of all points $(1,h,0)$ with 
$h{\hskip1pt\in\hskip1pt}\HH,\, h\overline h{=}1,\; hi{=-}ih$. Denote the polar 
lines of two points $p,p'{\,\in\,}U$  by $L$ and $L'$, and choose points 
$q{\,\in\,}pp'\smcap U$ and $z{\,\in\,}L\sm\cD$. 
As $\Delta_z\smcap\Kappa{\,=\,}\1$, we have 
$\Delta_z{\cong\,}\Delta/\Kappa$ , and representation of $\Delta_z$ 
on the compact group $\Kappa$ shows $\Delta{\,=\,}\Delta_z{\times}\Kappa$. 
The connected component $\Alpha$ of the center of $\Kappa$ is in the center of $\Delta$; it is not a Lie group.
The stabilizer $\Gamma{=\,}\Delta_{z,p,p',q}$ fixes also the point
$s{\,=\,}L\smcap L'{\,\in\,}\cD$, and $\dim\Gamma{\,=\,}4$.  
The fact that $\Delta_z$ is (doubly) transitive on $U$ implies that 
$\Delta_z$ has a subgroup $\Phi{\,\cong\,}\SU3\CC$ which is transitive on 
$U$, and $\Gamma'{\cong\,}\SU2\CC$. Let $\iota$ be the involution in
$\Gamma'$, obviously it is {\it planar\/}. It is a consequence of the stiffness 
result \cite{cp} 83.11 that the involution $\overline\iota{\,=\,}\iota|_\cD$ is not 
planar. Hence $\iota$ induces a reflection $\overline\iota$ on $\cD$.
Its axis $pp'\smcap\cD{\,\approx\,}\sS_4$ is also a line of $\cF{\,=\,}\cF_\iota$, 
and all lines of $\cF$ are homeomorphic to $\sS_4$. 
The central group $\Alpha$ acts effectively on $\cF$, and 
$\dim\Alpha{\,<\,}4$. By the structure theorem for compact groups 
(\cite{cp} 93.11), the semi-simple factor $\Kappa'$ is locally isomorphic to 
$(\SU2\CC)^2$. Exactly one of the two factors of $\Kappa'$ acts trivially 
on $\cF_\Gamma{\,=\,}\cF$, the other factor induces a $3$-dimensional Lie 
group on $L\smcap\cF$. By Richardson's theorem, $\Alpha$ would be a 
Lie group. \\
(f) Of the possibilities listed in step (b), only the following remains: $\cD$ is a Hughes plane (cf. \cite{cp} \S\,86, in particular, 86.33,34)\,), 
$\dim\Delta{\,=\,}17$, and there is an invariant subplane $\cC$ isomorphic to the classical complex plane. The commutator group 
$\Upsilon{=\,}\Delta'{\,\cong\,}\SL3\CC$ induces on $\cC$ the full automorphism group $\PSL3\CC$. The center of $\Upsilon$ is generated by an element 
$\epsilon$ of order $3$. The kernel $\Kappa$ is compact by stiffness, and 
$\dim\Kappa{\,\ge\,}5$.  Moreover, $\Upsilon$ acts trivially on $\Kappa$, and 
$\Kappa{\,\le\,}\Cs{}\Upsilon$. Choose a point $z$ which is not incident with 
any line of $\cC$. As $\cD$ is a Baer subplane, there is a unique line $L$ of 
$\cD$ with $z{\,\in\,}L$, and $\cF_\epsilon{\,=\,}\cC$ implies $L^\epsilon{\,\ne\,}L$.
The stabilizer $\Ypsilon_{\!z}$ fixes each point of $z^\Kappa{\,\subset\,}L$ 
and of $z^{\epsilon\Kappa}{\,\subset\,}L^\epsilon$. Both orbits of $z$ are homeomorphic to $\Kappa$; their dimension is at least $5$. Hence 
$\langle z^\Kappa,z^{\epsilon\Kappa}\rangle{\,=\,}\cP$ and 
$\Ypsilon_{\!z}{\,=\,}\1$. It follows that $z^{\!\Ypsilon}$ is open in $P$, and 
$\Delta$ would be a Lie group by 2.6. \Qed
\par\bigskip
{\Bf 4. Exactly one fixed element}
\par\medskip
Throughout this section, $\Delta$ fixes a line $W$ but no point.
\par\medskip
{\bf 4.0.} {\it If $\cF_\Delta{\,=\,}\{W\}$ and if $\dim\Delta{\,\ge\,}22$, then 
$\Delta$ is a Lie group\/}. 
\par\smallskip
{\tt Proof.} Accordimg to \cite{66} the assertion is true for 
$\dim\Delta{\ge\,}23$.\:Thus we may assume that 
$\dim\Delta{\,=\,}22$ and that 
$\Delta$ is not a Lie group. Steps (a,b) are the same as in \cite{66}. \\
(a) There exist arbitrarily small compact central  subgroups $\Nu{\,\le\,}\Delta$ of dimension $0$ such that $\Delta/\Nu$ is a Lie group, cf. \cite{cp} 93.8. 
If $\Nu$ acts freely on  $P\sm W$, then each stabilizer $\Delta_x$ with 
$x{\,\notin\,}W$  is a Lie group because $\Delta_x{\smcap}\Nu{\,=\,}\1$. 
By \cite{cp} 51.\hskip1pt6 and 8 and 52.12, the one-point compactification $X$ of $P\sm W$ is homeomorphic to the quotient space $P/W$, and $X$ is a Peano continuum (i.e., a continuous image of the unit interval); moreover, $X$ is homotopy equivalent to $\sS_{16}$, and  $X$ has Euler characteristic $\chi(X){\,=\,}2$.  
According to a theorem of L\"owen \cite{lwL}, these properties suffice for $\Delta$ to 
be a Lie group. \\
(b) Suppose now that  $x^\zeta{\,=\,}x$ for some $x{\,\notin\,}W$ and some 
$\zeta{\,\in\,}\Nu\sm\1$.  By assumption, $x^\Delta$ is not contained in a line and hence generates a $\Delta$-invariant subplane 
$\cD{\,=\,}\langle x^\Delta,W\rangle{\,\le\,}\cF_\zeta{\,<\,}\cP$. 
Put $\Delta^{\!*}{\,=\,}\Delta|_\cD$ and write 
$\Delta|_\cD{\,=\,}\Delta/\Kappa$. 
If $\cD$ is flat, then $\dim\Delta{\le\,}6{+}14$ by Stiffness; 
similarly,  $\dim\cD{\,=\,}4$ implies $\dim\Delta{\,\le\,}12{+}8$. 
Hence $\cD{\;\ldot}\cP$,\: $\Kappa$ is compact, $\dim\Kappa{\,<\,}8$,\: 
$\dim\Delta^{\!*}{\,\ge\,}15$,\: $\Delta^{\!*}$ is a Lie group\, (\cite{pwL}), 
and we may assume that $\1{\,\ne\,}\Nu{\,\le\,}\Kappa$. 
If even $\dim\Delta^{\!*}{\,>\,}16$, then \href{http://arxiv.org/abs/1402.0304}
{\cite{65} 1.10} shows that $\cD$ is the classical quaternion plane, because 
$\Delta^{\!*}$ does not fix a flag. \\
(c) If $\Delta$ is not transitive on $S{\,=\,}W\smcap\cD$, then there are
points $u,v,w{\,\in\,}S$ such that $\Delta{\hskip1pt:\hskip1pt}\Delta_{u,v,w}{\,\le\,}9$. Choose a line $L$ of $\cD$ in the pencil $\frL_v\sm\{W\}$ and a point 
$z{\,\in\,}L\sm\cD$, so that $\langle \cD,z\rangle{\,=\,}\cP$,\\  
$\Delta_z\smcap\Kappa{\,=\,}\1$, and $\Delta_z$ acts 
faithfully on $\cD$. In particular, $\Delta_z$ fixes $L$ and $\Delta_z$ is a Lie group. Moreover, $\dim z^{\Delta_L}{\,<\,}8$\: (or else $\Kappa$ would be a Lie group by  2.6 above.  
Therefore $\Delta{\hskip1pt:\hskip1pt}\Delta_z{\,\le\,}3{+}4{+}7{\,=\,}14$ and 
$\Lambda{\,=\,}(\Delta_{z,u,w})^1$ satisfies $\dim\Lambda{\,\ge\,}2$. 
The stiffness result 2.5(\^b) shows that $\Lambda{\,\cong\,}\Spin3\RR$.
From  \cite{cp} 55.32 and the structure of compact groups (\cite{cp} 93.11) it follows 
that  $\Kappa^1$ is commutative  or an almost direct product of 
$\Spin3\RR$ with a connected commutative group.
By  \cite{cp} 93.19 the group $\Lambda$ acts trivially on the commutative factor 
$\Alpha$ of $\Kappa^1$. As $\Kappa/\Nu\Kappa^1$ is finite and 
$\Kappa\smcap\Lambda{\,=\,}\1$, the group $\Lambda$ centralizes also a non-commutative factor of $\Kappa^1$, and hence $\Lambda{\,\le\,}\Cs{}\Kappa$. In the cases  $\dim\Kappa{\,\le\,}5$, $\dim\Delta_\cD{\,\ge\,}17$,  the plane $\cD$ is 
classical  and $\Lambda$ fixes all points of a circle in $S$ containing $u,v,w$.
If, on the other hand,  $\dim\Kappa{\,>\,}0$, then $z^\Kappa{\,\subseteq\,}\cF_\Lambda$ has positive dimension. In any case, $\cF_\Lambda$ is connected. 
Because $\Nu$ acts freely on $L\sm\cD$ and 
 groups of planes of dimension ${\le\,}4$ are Lie groups (cf.  \cite{cp} 32.21  and 71.2), it
 follows that $\cF_\Lambda{\,\ldot}\cP$. We have $\Kappa_z{\,=\,}\1$,\; 
 $\dim z^\Kappa{\,=\,}\dim\Kappa$,\; $z^\Kappa{\,\subseteq\,}\cF_\Lambda$, 
 and therefore $\dim\Kappa{\,\le\,}4$, in fact even $\dim\Kappa{\,<\,}4$ because of 2.6. Now 
$\dim\Delta^{\!*}{\,\ge\,}19$, and $\cD$ is the classical 
 quaternion plane, see \cite{cp} 84.28 or \cite{65} 1.10. In particular, 
$S{\,\approx\,}\sS_4$.    Recall that 
$\dim u^\Delta{\,+\,}\dim v^\Delta{\,+\,}\dim w^{\Delta_{u,v}}{\,\ge\,}8$.  
It follows that $\Delta$ is doubly transitive on the $3$-dimensional orbit 
$V{\,=\,}v^\Delta$. Similarly, $\Delta_z$ is doubly transitive on $V$, and 
$\dim w^{\Delta_{u,v}}{\,\ge\,}2$. The doubly transitive groups have been determined by Tits, see \cite{cp} 96.16. The large orbit of $w$ excludes the possibiliy that 
$V{\,\approx\,}\RR^3$. Hence $V$ is compact and 
$\Delta_z|_V$ is simple, in fact one of the groups $\PU2(\HH,1)$ or 
$\PSU3(\CC,1)$. In the first case, $\Lambda$ would be $4$-dimensional. 
In the second case, we obtain $\dim\Lambda{\,=\,}2$, contradicting stiffness. \\
(d) The previous step shows that $\Delta$ is transitive on $S{\,=\,}W\smcap\cD$. Hence 
$S{\,\approx\,}\sS_4$, and $\Delta$ has a subgroup {\cyss Yu}${\,\cong\,}\U2\HH{\,\cong\,}\Spin5\RR$, see  \cite{cp} 53.2,\;96.19--22,\;55.40, and 94.27. The central involution $\sigma$ of {\cyss Yu} induces on $\cD$ a reflection with axis $S$. As its center is not fixed by $\Delta$, it follows from \cite{cp} 61.13 that $\Delta^{\hskip-1.5pt*}$ contains a transitive translation group with axis~$S$. Therefore $\Delta{\hskip1pt:\hskip1pt}\Kappa{\,\ge\,}18$ and $\cD$ is classical  (cf. \cite{cp} p.\,500,\,Theorem or 
\cite{65} 1.10). By \cite{cp} 55.40 the group $\SO5\RR$ does not act on any compact plane. Hence the involution $\sigma$ is even a reflection of $\cP$, and 
$\sigma^\Delta\sigma{\,=\,}\Tau{\,\cong\,}\RR^8$ is normal in $\Delta$. 
The group  $\Delta|_\cD$ is contained in 
$\SL2\HH{\hskip1pt\cdot\hskip1pt}\HH^{\times}
{\ltimes}\RR^8$. 
Either  $\Delta|_S$ is compact, or $\Delta$ has a subgroup 
$\Ypsilon{\,\cong\,}\SL2\HH$\, (use \cite{cp} 94.34).
In the first case, $\Delta|_S{\,=\,}${\cyss Yu}$|_S{\,\cong\,}\SO5\RR$. For $v,w{\,\in\,}S$, 
the stabilizer $\Delta_v$ fixes a second point $u{\,\in\,}S$, and 
$\Delta{\hskip1pt:\hskip1pt}\Delta_{u,v,w}{\,=\,}4{+}3$. Choose  a line $L$ of $\cD$ in the pencil  $\frL_v\sm\{W\}$, a point $z{\,\in\,}L\sm\cD$, and put 
$\Gamma{\,=\,}\Delta_{L,w}$. Again we may assume that $\Nu{\,\le\,}\Kappa$. 
By 2.6 we have $\dim z^\Gamma{<\,}8$. The dimension formula yields 
$\Delta{\hskip1pt:\hskip1pt}\Gamma{\,\le\,}11$.
Let $\Lambda{\,=\,}(\Gamma_{\hskip-2pt z})^1$. Then $\dim\Lambda{\,\ge\,}4$.   
Again  $\Lambda$ is a Lie group since $\Nu$ acts freely on $L\sm\cD$. 
As $\cD$ is the classical quaternion plane, any collineation which fixes $3$ collinear points of $\cD$ even fixes all points of a circle. Hence $\cF_\Lambda$ is connected. 
By \cite{cp} 32.21  and 71.2, we  have $z^\Nu{\,\subseteq\,}\cF_\Lambda{\;\ldot}\cP$. 
Now Stiffness would imply $\dim\Lambda{\,\le\,}3$ contrary to what has been stated before. \\
(e) If $\Delta$ has a subgroup $\Ypsilon{\,\cong\,}\SL2\HH$, then the arguments in 
step~(d) imply $\Delta{:}\Kappa{\,\ge\,}15{+}8$. A maximal compact subgroup 
of $\Delta$ is connected and properly contains {\cyss Yu}. Consequently 
$\dim\Delta{\,\ge\,}24$, and $\Delta$ is a Lie group by \cite{66} 4.0. \Qed
\par\medskip
For $8$-dimensional planes the following has been proved in 
\cite{65} 3.1:
\par\smallskip
{\bold 4.1.}  {\it Suppose that $\cF_\Delta$ consists of a unique line in an \emph{$8$-dimensional plane} and that $\Delta$ is a semi-simple group of dimension $>\!3$.  Then $\Delta$ is a Lie group and $\dim\Delta{\,\le\,}10$. In the case of equality, $\Delta$ is isomorphic to  $\Opr5(\RR,1)$ or to some covering group of $\Opr5(\RR,2)$\/}.
\par\medskip
{\bold 4.2.} 
{\it Let $\cF_\Delta{\,=\,}\{W\}$. If $\Delta$ is semi-simple and if $\dim\Delta{\,\ge\,}11$, then $\Delta$ is a Lie group\/}.
\par\smallskip
{\tt Proof.} As in steps (a,b) of the previous proof we may consider a subplane 
$\cD{\,=\,}\langle x^\Delta,W\rangle<\cP$. Both groups
$\Delta^{\!*}{\,=\,}\Delta|_\cD{\,=\,}\Delta/\Kappa$ and 
$\Lambda{\,=\,}\Kappa^1$ are semi-simple. 
If $\cD$ is flat and $\Delta^{\!*}{\,\ne\,}\1$, then 
$\Delta^{\!*}{\,\circeq\,}\SL2\RR$ and 
 either contains a central reflection or acts as hyperbolic motion group without fixed element. Hence 
$\Delta^{\!*}{\,=\,}\1$ and 
$\Delta{\,=\,}\Lambda{\,=\,}\Gtwo$ or 
$\dim\Delta{\,\le\,}10$ by stiffness. 
If $\dim\cD{\,=\,}4$, then \cite{66} 2.14 shows that 
$\dim\Delta^{\!*}{\,\le\,}3$ and $\Lambda{\,\cong\,}\SU3\CC$ 
or $\dim\Lambda{\,<\,}8$ by stiffness. Thus $\Delta$ is
a Lie group by 2.9 or $\dim\Delta{\,\le\,}10$. In the case 
$\cD{\,\,\ldot}\cP$ the kernel $\Kappa$ is compact 
and then $\Lambda{\,\cong\,}\Spin3\RR$ and 
$\dim\Delta^{\!*}{\ge\,}8$.  By the result 4.1 the group 
$\Delta^{\!*}$ is a Lie group and 
$\dim\Delta^{\!*}{\le\,}10$. From 2.9 it follows that $\Delta$  is a Lie group. \Qed
\par\bigskip
{\Bf 5. Fixed flag} 
\par\medskip
For $8$-\emph{dimensional} planes, the following results have been proved in \cite{63} Theorem 1.3 and in \cite{69} 5.1:
\par\medskip
{\bf 5.0.} {\it If the semi-simple group $\Delta$ of an $8$-\emph{dimen\-sional} plane fixes exactly one line and possibly some points on this line, and if $\dim\Delta{\,>\,}3$, then $\Delta$ is a Lie group\/}.
\par\medskip
{\bf 5.1.} {\it Let $\cF_\Delta{\,=\,}\{v,W\}$ be a flag in an $8$-\emph{dimen\-sional} plane. If $\dim\Delta{\,\ge\,}10$, then $\Delta$ is a Lie group\/}.
\par\medskip
{\bf 5.2.} {\it If $\Delta$ fixes a flag 
$\{v,W\}$  and possibly one or two further points on $W$, but no other line, if $\Delta$ is semi-simple,  and if  
$\dim\Delta{\,\ge\,}11$, then $\Delta$ is a  Lie group\/}.
\par\smallskip
{\tt Proof.} The first arguments are the same as in steps (a,b) of the proof of 4.0 above. Thus there is a connected proper subplane 
$\cD{\,=\,}\langle x^\Delta,v,W\rangle$ of $\cP$. 
Put $\Delta^{\!*}{\,=\,}\Delta|_\cD{\,=\,}\Delta/\Kappa$, and note that both 
$\Delta^{\!*}$ and $\Lambda{\,=\,}\Kappa^1$ are semi-simple. If $\cD$ is flat, then $\Delta^{\!*}$ is trivial by \cite{cp} 33.8, the kernel $\Kappa$ is isomorphic to the Lie group 
$\Gtwo$ or $\dim\Kappa{\,\le\,}10$ by stiffness 2.5(d). 
If $\dim\cD{\,=\,}4$, then $\dim\Delta^{\!*}{\,\le\,}3$\; (see \cite{66} 2.14 or \cite{cp} 71.\,8,10 and 72,\,1--4), and stiffness 2.5(a,e) shows that 
$\dim\Lambda{\,<\,}8$ or $\Lambda{\,\cong\,}\SU3\CC$ and 2.9 applies. \\ 
Finally, let  $\cD{\,\,\ldot}\cP$.  Then  
$\dim\Delta^{\!*}{\le\,}3$ or $\Delta^{\!*}$ is a Lie group by 
5.0 above. Moreover, $\Lambda$ is compact and semi-simple, hence a Lie group, and then 
$\Lambda{\,\cong\,}\Spin3\RR$ by 2.5(\^b) and 
$\dim\Delta^{\!*}{>\,}3$. Therefore $\Delta$  is a Lie group by 2.9. \Qed
\par\medskip
\par\bigskip
{\Bf 6. Two non-incident fixed elements}
\par\medskip
It has been stated in 2.3 above that $\Delta$ is a Lie group whenever 
$\dim\Delta{\,\ge\,}27$. In the case $\cF_\Delta{\,=\,}\{a,W\}, \ a{\,\notin\,}W$,
no improvement of this result has been found in general; for semi-simple groups  
it has been proved in \cite{psz}\hskip1pt({\bold a}) that  
$\dim\Delta{\,\ge\,}26$ suffices.  \par
Suppose that $\Zeta{\,=\,}\Cs{}\Delta$ is not a Lie group. If $\Zeta|_W{\,\ne\,}\1$, then there is a point $p$ such that $\Lambda{\,=\,}(\Delta_p)^1$ fixes a 
quadrangle, and stiffness implies $\dim\Lambda{\,\le\,}11$. On the other hand
 $\Delta{:}\Delta_p{\,<\,}16$ by 2.6. Thus 
 $\dim\Delta{\,\le\,}26$.
 \par\medskip
Throughout  section $6$, assume that $\Delta$ fixes a 
 line $W$ and a point $a{\,\notin\,}W$.
\par\medskip
{\bold 6.1 Proposition.}  {\it If $\Delta$  is semi-simple, and if there is  some \emph{planar} element $\zeta{\,\in\,}\Zeta$ {\rm(}so that $\cF_\zeta$ is a  subplane 
of $\cP${\rm)}, then $\dim\Delta{\,\le\,}19$\/}.
\par\smallskip
{\tt Proof.}  Assume that $\dim\Delta{\,>\,}19$, and let
$x^\zeta{\,=\,} x{\,\in\,}\cF_\zeta$. Then $x^\Delta{\,\ne\,}x$ 
by assumption, $x^\Delta{\,\subseteq\,}\cF_\zeta$, and $\cF_\zeta$ is connected.    
Note that $\Delta|_{\cF_\zeta}{\,=\,}\Delta/\Kappa$ and 
$\Kappa^1$ are semi-simple. If $\cF_\zeta$ is flat or if 
$\dim\cF_\zeta{\,=\,}4$, then stiffness yields $\dim\Kappa{\,\le\,}14$, and
$\Delta{:}\Kappa{\,\le\,}3$\; (see \cite{66} 2.14 or \cite{cp} 71.8--72.4). 
If $\cF_\zeta{\,\,\ldot}\cP$, then $\Kappa^1$ is compact and semi-simple, hence a Lie group, and $\Kappa^1{\,\cong\,}\Spin3\RR$ by 2.5(\^b).
 The classification of all planes such that $\Delta{:}\Kappa{\,>\,}16$
 as summarized in \cite{65} 1.10 shows that $\cF_\zeta$ is the classical quaternion plane, and then $\Delta/\Kappa{\,\cong\,}\SL2\HH{\cdot}\Spin3\RR$.  It follows that a maximal compact subgroup of $\Delta$ is locally isomorphic to
$\Spin5\RR{\cdot}(\Spin3\RR)^2$, and $\Delta$ would be a Lie group. 
Thus $\dim\Delta{\,\le\,}16{+}3$. \Qed
\par\medskip
{\bold 6.2 Proposition.}  {\it Suppose that $\Delta$  is semi-simple and that $\dim\Delta{\,\ge\,}25$. If the center $\Zeta$  of $\Delta$ acts non-trivially on $W$, then $\Delta$ is a Lie group\/}.
\par\smallskip
{\tt Proof.} Because of  \cite{psz}\hskip1pt({\bold a}) we may assume that $\dim\Delta{\,=\,}25$. Assume that $\Zeta$ is not a Lie group.  There are several possibilities: \\
(a) $\Delta$ has an almost simple factor $\Psi$ of dimension $16$. Then $\Psi$ is locally isomorphic to $\SL3\CC$ and 
$\Psi$ is a Lie group. Denote the product of the other factors of $\Delta$ by $\Gamma$.
Let $v^\Zeta{\,\ne\,}v{\,\in\,}W$ and  
$p{\,\in\,}K{:=\,}av\sm\{a,v\}$.  Then
$\Lambda{\,=\,}(\Psi_p)^1$ has positive dimension by 2.6,
and $\cD{\,=\,}\cF_\Lambda$ is a proper subplane. 
In particular, $\dim\cD{\,\in\,}\{0,2,4,8\}$. 
Put $\Xi{\,=\,}(\Gamma_{\!p})^1$ and note that 
$p^\Gamma$ is contained in $\cD$. It follows that 
$\dim\Xi{\,>\,}0$. Hence $\cE{\,=\,}\cF_\Xi$ is also a proper subplane. Now $\cE^\Psi{\,=\,}\cE$ and 
$\dim\Lambda{\,\ge\,}8$. Stiffness implies 
$\dim\cD{\,\le\,}4$. If  $\cD$ is connected, then $\Zeta$ 
induces a Lie group on $\cD$ and there is some 
$\delta{\,\in\,}\Zeta$ such that 
$\cD{\,\le\,}\cF_\delta{\,=\,}\cF{\,=\,}\cF^\Delta{\,<\,}\cP$. 
The factor $\Psi$ acts almost effectively on $\cF$. 
Hence $\cF{\,\ldot}\cP$\, (use \cite{cp} 71.8 and 72.8).
Let $\Delta|_\cF{\,=\,}\Delta/\Kappa$. As $\Kappa^1$	 is 
semi-simple, stiffness shows $\dim\Kappa{\,\le\,}6$ and
$\Delta{:}\Kappa{\,=\,}19$. By the classification of
$8$-dimen\-sional planes as summarized in 
\cite{65} 1.10 (cf. also \cite{cp} 84.28), the plane $\cF$ is the classical quaternion plane and 
$\Delta|_\cF{\,\cong\,}\SL2\HH{\cdot}\HH^{\times}$ in 
contradiction to the action of $\Psi$ on $\cF$. 
Consequently $\dim\cD{\,=\,}0$,\, 
$\dim\Xi{\,=\,}9$, and $\dim\cE{\,\le\,}4$ by stiffness, 
but then $\Psi|_\cE{\,=\,}\1$, which is impossible.   \\
(b) $\Delta$ has an almost simple factor $\Psi$ of dimension $15$, and $\Psi$ is locally isomorphic to one of the groups $\SU4\CC$, $\SL2\HH$, $\SL4\RR$, or 
$\Psi$ maps onto a group $\PSU4(\CC,r)$, 
$r{\,=\,}1,2$. Then there is a second almost simple factor 
$\Gamma$ of  dimension $10$. 
This situation can be treated in a similar way as case (a), but a few more arguments are needed. \\
($i$)  {\it $\Delta$ is not  transitive on $W$\/}. 
In fact, trans\-itivity implies $W{\,\approx\,}\sS_8$ by 2.6,  
a maximal  compact subgroup $\Phi$ of $\Delta|_W$ is also transitive (\cite{cp} 96,19), and  $\Phi{\,\cong\,}\SO9\RR$ would be too big  (cf. \cite{cp} 96.23). \\
($ii$) {\it There are points $v{\,\in\,}W$ and 
$p{\,\in\,}K{\,:=\,}av\sm\{a,v\}$ such that 
$\dim v^\Delta{\,<\,}8$ and $\Psi_p{\,\ne\,}\1$\/}. \\
The first claim follows from ($i$), and 
$\dim\Psi_v{\,\ge\,}8$. If $\Psi_v$ acts freely on 
$K$, then $\Psi_v$ is transitive 
on $K$ by \cite{cp} 96.11, lines are homeomorphic to 
$\sS_8$ by 2.6, and 
$\Psi_v{\,\approx\,}K{\,\simeq\,}\sS_7$ (see \cite{cp} 52.5). 
A maximal compact subgroup $\Phi$ of $\Psi_v$ is at most $7$-dimensional and has the homotopy of $\sS_7$,  but 
$\pi_\nu\sS_7{\,=\,}0$ for $\nu{\,<\,}7$ and 
$\pi_\nu\Phi{\,\ne\,}0$ for $\nu{\,=\,}{\rm1\ or\ 3}$,  see \cite{cp}.  \\
($iii$) {\it One of $\cD{\,=\,}\cF_{\Psi_p}$ or 
$\cE{\,=\,}\cF_{\Gamma_{\!p}}$ is a connected proper subplane\/}. Note that $p^\Gamma{\subseteq\,}\cD$. 
If $v^\Gamma{\,\ne\,}v$, then $\cD$ contains a quadrangle 
and a connected orbit of $\Gamma$; thus $\cD{\,<\,}\cP$.
If $v^\Gamma{\,=\,}v$, however, then 
$v^\Psi{\,\ne\,}v$, and $\Gamma$ fixes each point of the 
connected orbit $v^\Psi{\,\subseteq\,}W$; moreover, 
$\dim\Gamma_{\!p}{\,>\,}1$, and $\cE{\,<\,}\cP$. \\
($iv$) If $v^\Gamma{\ne\,}v$, then 
$\dim\Gamma|_\cD{\,=\,}10$,\,   $\cD{\,\ldot}\cP$, 
and $\dim\Gamma_{\!p}{\ge\,}2$.  We have 
$\cE^\Psi{\,=\,}\cE$. Either $\cE{\,<\,}\cP$ or 
$v^\Psi{\,=\,}v$.  In the latter case, $\Psi_p|_\cD{\,=\,}\1$,\, 
$\dim\Psi_p{\,=\,}7$ by stiffness, and 
$\Psi$ is transitive on $K$. Hence 
$\Gamma_{\!p}{\,\le\,}\Delta_{[av]}$, and 
$\Gamma_{\!p}{\,\triangleleft\,}\Gamma$ would be a proper 
normal subgroup. Consequently 
$\dim\Psi|_\cE{\,=\,}15$, and $\cE{\,\ldot}\cP$. \\
($v$) If $v^\Gamma{=\,}v$, then 
$v^\Psi{\,\ne\,}v$,\, $\cE{\,<\,}\cP$,\, 
$\dim\Psi|_\cE{\,=\,}15$, and $\cE{\,\ldot}\cP$. 
From \cite{63} 3.3 and $\dim\Psi{\,>\,}10$ it follows that 
$\cE$ is the classical quaternion plane. In particular, 
$\Psi{\,\cong\,}\SL2\HH$. Moreover 
$\Zeta|_\cE$ is a Lie group and there is some 
$\zeta{\,\in\,}\Zeta$ such that 
$\cE{\,=\,}\cF_\zeta$ is $\Delta$-invariant and 
$\dim\Delta|_\cE{\,=\,}25$. This is impossible. \\
($vi$) {\tt Corollary.} {\it Both $\cD$ and $\cE$ are Baer 
subplanes\/}.  \\
Again $\cE$ is classical by 2.11. As in the previous step, $\cE^\Delta{\,=\,}\cE$ and $\dim\Delta|_\cE{\,=\,}25$, 
a final contradiction. \\
(c) $\Delta$ has an almost simple factor $\Psi$ of dimension $14$. In this case, $\Psi$ is a Lie group, there are two other factors  $\Beta$ and $\Alpha$ of dimension $8$ and $3$, 
respectively, and $\Beta$ is a a Lie group or maps onto 
$\PSU3(\CC,1)$.  \\ 
($i$) If $\Psi{\,\cong\,}\Gtwo$ is compact or, more generally,  whenever $\Psi$ has a subgroup 
$\Phi{\,\cong\,}\SO3\RR$, then any two conjugate 
commuting  involutions $\iota,\kappa{\,\in\,}\Phi$ are
planar (use \cite{cp} 55.35), and $\cF_\iota$ is a Baer subplane. (For involutions in $\Gtwo$ see \cite{cp} 31-35 or \cite{66} 2.5.)  Choose a point 
$x$ in $\cF_\iota,\ a{\,\ne\,}x{\,\notin\,}W$, and put 
$v{\,=\,}ax\smcap W$ and $\Gamma{\,=\,}\Beta\Alpha$. 
Either $\cF_{\iota,\kappa}$ is a $4$-dimensional subplane or 
$\kappa$ induces on $\cF_\iota$ a  reflection with center 
$v$, say. In the first case, $\Zeta|_{\cF_{\iota,\kappa}}$ is 
a Lie group by \cite{cp} 71.2, and there is some 
$\zeta{\,\in\,}\Zeta$ such that $\cF_\zeta$ is a connected 
proper subplane, but then $\Delta|_{\cF_\zeta}$ 
is too big. In the second case, $v^\Gamma{\,=\,}v$,\,   
$V{\,=\,}v^\Psi{\,\ne\,}v$, and $\Gamma|_V{\,=\,}\1$. 
It follows that $\dim\Gamma_{\!x}{\ge\,}3$ and  that
$\cE{\,=\,}\cF_{\Gamma_{\!x}}{\,=\,}\cE^\Psi$ is a connected 
proper subplane. The simple group $\Psi$ acts effectively on $\cE$. Consequently $\iota|_\cE$ is planar and 
$\dim(\cE\smcap\cF_\iota){\,=\,}4$. Again this is impossible, 
cf. also \cite{cp} 84.8. \\
($ii$) Similar arguments apply if $\Beta$ has a subgroup 
$\Phi{\,\cong\,}\SO3\RR$, in particular, if $\Beta$ is compact. Again involutions $\iota,\kappa{\,\in\,}\Phi$ are planar and $\cF_{\iota,\kappa}$ is a triangle. The factor
$\Psi$ fixes all points of some $\Beta$-orbit  
$V{\,=\,}v^\Beta{\,\subseteq\,}W$, and $\dim\Psi_x{\,\ge\,}6$. 
The center $\Zeta$ acts effectively on 
$\cD{\,=\,}\cF_{\Psi_x}$, and $\cD{\,\ldot\,}\cP$. Now 
stiffness 2.5(\^b) shows $\dim\Psi_x{\,\le\,}3$ because 
$\Psi$ is  Lie group,  a contradiction. Therefore $\Beta$ has no subgroup
$\SO3\RR$. \\  
($iii$) Therefore $\Psi$ is isomorphic to the twofold 
coverring $\tilde{\rm G}$ of $\Gtwo(2)$; a maximal compact subgroup $\Phi$ of $\Psi$ is   isomorphic to $\Spin4\RR$.
The center of $\Psi$ is generated by a reflection 
$\omega$ with axis $W$. There are two other central involutions  $\iota,\kappa$ in $\Phi$, these may be 
planar or they are reflections with centers $u,v{\,\in\,}W$. In the first case the argumentation is the same as in step ($i$).
If $\iota{\,\in\,}\Delta_{[u,av]}$, then 
$v^\Gamma{\,=\,}v{\,\ne\,}v^\Psi{\,=\,}V$ and 
$\Gamma|_V{\,=\,}\1$. For the same reasons as in 
step (b)($ii$) the group $\Beta$ does not act freely on
$K{\,=\,}av\sm\{a,v\}$ and $\Beta_x{\,\ne\,}\1$ for a suitable 
$x{\,\in\,}K$. It follows that   
$\cE{\,=\,}\cF_{\Beta_{\!x}}{\,=\,}\cE^{\Psi\Alpha}$ is  a proper connected sub plane. Because 
$\omega$ is a reflection of $\cP$, the group $\Psi$ acts 
effectively on $\cE$. The arguments can now be repeated for $\Psi|_\cE$. Hence $\iota$ and $\kappa$ induce reflections on $\cE$, and $\Gamma$ fixes all points of 
$V{\,\subseteq\,}W$. Moreover, \cite{cp} 96.13, applied to 
$\Phi|_V$, shows that $\dim V{\,\ge\,}3$, so that 
$\cE{\,\ldot}\cP$. If $\dim(\Psi\Alpha)|_\cE{\,=\,}17$, 
the classification as summarized in
 \cite{cp} 84.28 or \cite{65} 1.10 would imply that $\cE$ 
 is a translation plane. Therefore $\Alpha|_\cE{\,=\,}\1$ 
 and $\cE{\,=\,}\cF_\Alpha$ is $\Psi\Beta$-invariant.
By stiffness $\dim(\Psi\Beta)|_\cE{\,=\,}22$ is too big. 
Thus case (c) is impossible. \\
(d) An almost simple factor $\Psi{\,\triangleleft\,}\Delta$ of maximal dimension satisfies $\dim\Psi{\,=\,}10$; all 
other factors have dimension $3$ or $6$. 
Write $\Delta$ as an almost direct product 
$\Beta\Alpha$ with $\dim\Alpha{\,=\,}3$, and choose
$v{\,\in\,}W$ such that   $v^\Zeta{\,\ne\,}v$. Let 
$a,v{\,\ne\,}x{\,\in\,}av$. Then 
$\dim\Beta_x{>\,}6$\,  (by 2.6), and 
$\cE{\,=\,}\cF_{\Beta_x}{\,=\,}\cE^\Alpha$  
is a proper subplane.  If $\dim\cE$ is $2$ or $4$, then 
$\Zeta|_\cE$ is a Lie group, and there is some 
$\zeta{\,\in\,}\Zeta$ such that
$\cE{\,\le\,}\cF{\,=\,}\cF_{\hskip-2pt\zeta}{\,=\,}
\cF^\Delta{\,\ldot}\cP$. By stiffness  
$\dim\Delta|_\cF{\,=\,}19$,\, $\cF$ is the classical quaternion plane (cf. \cite{cp} 84.28 or \cite{65} 1.10), and then
$\Delta|_\cF{\,\cong\,}\SL2\HH{\cdot}\HH^{\times}$ is not 
semi-simple. Therefore $\dim\cE{\,\in\,}\{0,8\}$.  
If $\cE$ is not connected, then $\Alpha|_\cE{\,=\,}\1$ and 
$\cD{\,=\,}\cF_\Alpha{\,=\,}\cD^\Beta{\,=\,}\cD^\Delta$. 
By stiffness, the semi-simple kernel  of 
$\Delta|_\cD$ has dimension at most $6$  and  contains  $\Alpha$. Hence $\dim\Delta|_\cD{\,\ge\,}19$, but   
again this is impossible. Thus  
$\cE{\,=\,}\cE^\Zeta{\,\ldot}\cP$, and $\Zeta$ acts 
effectively on $\cE$\, (or else $\cE{\,=\,}\cF_\zeta$ for some 
$\zeta{\,\in\,}\Zeta$,\, $\cE^\Delta{\,=\,}\cE$,\, 
$\dim\Delta|_\cE{\,=\,}19$, a contradiction.)
Now $\Beta_x\smcap\Zeta{\,=\,}\1$,\, 
$\Beta_x$ is a Lie group, and $\dim\Beta_x{\,\le\,}3$ by stiffness, in contrast to what has been stated above. \\
(e) $\Delta$ has exactly two factors of dimension $8$; all
other factors have dimension $3$ or $6$. This is the last possibility, it can be treated exactly in the same way 
as (d). \\
{\tt Remark.} The cases (a) and (c) can also be exclude 
as in step (d). \Qed
\par\medskip
{\bold 6.3 Proposition.}  {\it Suppose that $\Delta$  is semi-simple and that $\dim\Delta{\,\ge\,}25$. If the center $\Zeta$  of $\Delta$ acts trivially on $W$, then $\Delta$ is a Lie group\/}. 
\par\smallskip
{\tt Proof.} Again we may assume that $\dim\Delta{\,=\,}25$. Suppose that $\Zeta$ is not a Lie group.  There are the same possibilities as in 6.2. \\
(a) $\Delta$ has an almost simple factor $\Psi$ of dimension $16$. Then $\Psi$ is locally isomorphic to $\SL3\CC$ and 
$\Psi$ is a Lie group. Denote the product of the other factors of $\Delta$ by $\Gamma$. The factor $\Psi$ has a 
subgroup $\Phi{\,\cong\,}\SO3\RR$. As before, the involutions in $\Phi$ are necessarily planer. Let 
$\iota,\kappa$ be two commuting involutions in $\Phi$. 
If $\kappa$ induces a planar involution on $\cF_\iota$, then 
$\Zeta$  acts freely on the $4$-dimensional plane 
$\cF_{\iota,\kappa}$, and $\Zeta$ would be a Lie group. 
Therefore $\kappa|_{\cF_\iota}$ is a reflection, and 
$\kappa$ fixes points $u,v{\,\in\,}W\smcap\cF_\iota$. 
Choose $x{\,\in\,}\cF_\iota$ not incident with the three fixed lines. Then $\dim x^\Gamma{<\,}8$ because of 2.6, and
$\dim\Gamma_{\!x}{\ge\,}2$. The group $\Gamma_{\!x}$ fixes a quadrangle and 
$\cE{\,=\,}\cF_{\Gamma_{\!x}}{\,=\,}\cE^\Psi$ is a proper subplane. We have $\dim\Psi|_\cE{\,=\,}16$ and 
$\cE{\,\ldot}\cP$. Hence $\cE$ is the classical quaternion 
plane by 2.11 . But then 
$\SL3\CC$ does not act on $\cE$ and the group $\Zeta$ 
would be a Lie group. \\
{\tt Remark.} The proof of (a) shows that a factor $\Psi$ 
of dimension ${>\,}10$ does not have a subgroup $\SO3\CC$. \\
(b)  $\Delta$ is an almost direct product $\Psi\Gamma$ of two factors  of dimension $15$ and $10$, respectively. 
As in 6.2(b)($i$,$ii$), $\Delta$ is not transitive on $W$, 
and there is a point $x$ on a line $av$ such that 
$\Psi_x{\,\ne\,}\1$. \\
If $V{\,=\,}v^\Gamma{\,\ne\,}v$, then 
$\cD{\,=\,}\cF_{\Psi_x}{=\,}\cD^\Gamma$ is a connected 
proper subplane, $\Zeta$ acts freely on $\cD$, and 
$\cD{\,\ldot}\cP$. Because of 2.6 it    
follows that $\dim x^\Gamma{<\,}8$ 
and  $\dim\Gamma_{\!x}{\ge\,}3$. 
As $\Gamma$ acts almost effectively on $V$ and 
$\Gamma$ has a compact subgroup of dimension at least $4$, Mann's theorem \cite{cp} 96.13 shows that
$\dim V{\,\ge\,}3$. Suppose that  $v^\Psi{=\,}v$. Then 
$\Psi_x|_V{\,=\,}\1$, $\dim x^\Psi{\,<\,}8$ by 2.6, and 
$\dim\Psi_x{\ge\,}8$, but $\cF_{\Psi_x}{\,\ldot}\cP$, 
and $\dim\Psi_x{\,\le\,}7$ by stiffness.Therefore 
$v^\Psi{\ne\,}v$ and 
$\cE{\,=\,}\cF_{\Gamma_{\!x}}{=\,}\cE^\Psi{\,\ldot}\cP$.  
In fact, $\cE$ is classical by 2.11. As a group of homologies $\Zeta$ acts freely on $\cE$, and $\Zeta$ would be a Lie group. Consequently $v^\Gamma{=\,}v$ whenever 
$\dim v^\Delta{\,<\,}8$. Moreover
$\dim\Gamma_{\!x}{\,>\,}1$ and  
$\cE{\,=\,}\cE^\Psi{\,\ldot}\cP$. Again $\cE$ is classical 
by 2.11, but then $\Zeta$ would be 
a Lie group, a contradiction. \\
(c) $\Delta$ has an almost simple factor $\Psi$ of dimension $14$. In this case, $\Psi$ is a Lie group, there are two other factors  $\Beta$ and $\Alpha$ of dimension $8$ and $3$, 
respectively, and $\Beta$ is a a Lie group or maps onto 
$\PSU3(\CC,1)$. This case is similar to the previous one.  
Put $\Upsilon{\,=\,}\Psi\Alpha$. Again $\Delta$ is not 
transitive on $W$, and $\dim v^\Delta{\,<\,}8$ for somme 
$v{\,\in\,}W$. Let $a,v{\,\ne\,}x{\,\in\,}av$. Then 
$\dim\Upsilon_{\!x}{>\,}1$. If $v^\Beta{\,\ne\,}v$, then    
$\cD{\,=\,}\cF_{\Upsilon_{\!x}}{=\,}\cD^\Beta{\,\ldot}\cP$, 
and $\dim\Beta_x{>\,}0$ because of 2.6. It follows  
that $v^\Upsilon{\ne\,}v$\, (or else 
$\dim\Upsilon_{\!x}{>\,}8$ and 
$\Upsilon_{\!x}|_\cD{\,=\,}\1$, which contradicts  stiffness). Hence
$\cE{\,=\,}\cF_{\Beta_x}{=\,}\cE^\Upsilon{\,\ldot}\cP$. 
The kernel of $\Upsilon|_\cE$ has dimension at most $3$,
and $\dim\Psi|_\cE{\,=\,}14$. Now  $\cE$ is classical 
by 2.11, and  $\Zeta$ would be a Lie group. 
If $v^\Beta{\,=\,}v$, however,  then 
$\Beta|_{v^\Delta}{\,=\,}\1$, and 
$\dim x^\Beta{\,<\,}8$\, (use 2.6). Again 
$\dim\Beta_x{>\,}0$, and $\cE$ is the classical quaternion plane by 2.11. As before, $\Zeta$ is a Lie group 
in contrast to the assumption. \\
(d) An almost simple factor $\Psi{\,\triangleleft\,}\Delta$ of maximal dimension satisfies $\dim\Psi{\,=\,}10$; all 
other factors have dimension $3$ or $6$. Write $\Delta$ 
as an almost direct product $\Upsilon\Gamma$ with 
$\dim\Upsilon{\,=\,}16$. The arguments in  6.2(b)($i$,$ii$) 
show that $\Delta$ is not transitive on $W$. Choose $v$ 
and $x$ as in step (c). By 2.6 we conclude that 
$\dim\Upsilon_{\!x}{>\,}1$. If $v^\Gamma{\,\ne\,}v$, then 
$\cD{\,=\,}\cF_{\Upsilon_{\!x}}{=\,}\cD^\Gamma{\,\ldot}\cP$ 
and $\Upsilon_{\!x}|_\cD{\,=\,}\1$. Stiffness implies 
$\dim\Upsilon_{\!x}{\le\,}7$ and $\dim x^\Upsilon{\ge\,}9$.
As $x^\Gamma$ is contained in $\cD$, we have 
$\dim\Gamma_{\!x}{\,>\,}0$. Therefore 
$\cE{\,=\,}\cF_{\Gamma_{\!x}}{\,=\,}\cE^\Upsilon{\ldot}\cP$ and $\dim x^\Upsilon{\le\,}8$. This contradiction shows that 
$v^\Gamma{=\,}v$ and $\Gamma|_{v^\Delta}{\,=\,}\1$.
Consequently again
$\cE{\,=\,}\cF_{\Gamma_{\!x}}{\,=\,}\cE^\Upsilon{\ldot}\cP$.
Let $\Upsilon|_\cE{\,=\,}\Upsilon/\Kappa$. Then 
$\Lambda{\,=\,}\Kappa^1$ is compact and semi-simple, 
hence a Lie group, and stiffness 2.5(\^b) shows\, 
$\dim\Lambda{\,\le\,}3$. As $\Zeta$ acts effectively on 
$\cE$, the plane $\cE$ is not classical. We have 
$\dim\Upsilon|_\cE{\,\ge\,}13$, and 2.11  implies that 
$\Upsilon/\Kappa{\,\cong\,}\U2(\HH,r){\cdot}\SU2\CC$. The complement of $\Psi$ in $\Upsilon$ can be chosen arbitrarily among the factors of dimension $3$ or $6$. 
Hence all these factors are compact and $3$-dimensional.
Therefore $\Delta$ is a Lie group. \\
(e) $\Delta$ has exactly two factors of dimension $8$; all
other factors have dimension $3$ or $6$. The first part of the proof is identical to the previous one, and there is a Baer subplane $\cE{\,=\,}\cF_{\Gamma_{\!x}}{\,=\,}\cE^\Upsilon$. 
As $\Upsilon$ has two $8$-dimensional factors, 
$\dim\Upsilon|_\cE{\,=\,}16$, and 2.11 implies that $\cE$ is the classical quaternion plane, but then $\Zeta$ is a Lie group contrary to the assumption. \Qed
\par\medskip
Combined, Propositions 6.2 and 6.3 yield 
\par\smallskip
{\bold 6.4 Theorem.} 
{\it If $\cF_\Delta{\,=\,}\{a,W\},\,a{\notin}W$, if $\Delta$ is 
semi-simple, and if $\dim\Delta{\,\ge\,}25$, then $\Delta$ is a Lie group\/}.
\par\bigskip
{\Bf 7. Fixed double flag}
\par\medskip
{\bold 7.1.} 
{\it Let $\cF_\Delta$ be a double flag $\{u,v,uv,ov\}$. If $\Delta$ is semi-simple and if $\dim\Delta{\,\ge\,}25$, then 
$\Delta$ is a Lie group\/}.
\par\smallskip
{\tt Proof.} For $\dim\Delta{\,\ge\,}26$, this has been 
proved in \cite{66} 8.0. It suffices, therefore, to consider the 
case $\dim\Delta{\,=\,}25$. The arguments are simi\-lar to 
those in \cite{66}. Let $\Psi$ be an almost simple factor of $\Delta$ of \emph{maximal} dimension. Then 
$\dim\Psi{\,\in\,}\{16,15,14,10,8\}$. The product of the other factors of $\Delta$  will be denoted by $\Gamma$, the center of $\Delta$ by $\Zeta$.  \\
(a) {\it If $p{\,\notin\,}uv\smcup ov$, then $\Zeta_p{\,=\,}\1$\/}: suppose that $p^\zeta{=\,}p$ for some 
$\zeta{\,\in\,}\Zeta\sm\{\1\}$. 
Then $\zeta|_{p^\Delta}{\,=\,}\1$ 
and $\cD{\,=\,}\langle p^\Delta\rangle{\,\ldot}\cP$. 
With $\Delta|_\cD{\,=\,}\Delta/\Kappa$ we get 
$\Delta{:}\Kappa{\,\le\,}$ and $\dim\Kappa{\,\le\,}7$, 
a contradiction. \\
(b) {\it Each reflection $\sigma{\,\in\,}\Psi$ 
is in $\Psi_{[u,ov]}$, and each non-central involution 
$\beta{\,\in\,}\Psi$ is planar\/}. In fact, 
$\sigma{\,\in\,}\Psi_{[uv]}$ would imply that 
$\sigma^\Delta\sigma{\,\cong\,}\RR^k$ is normal in 
$\Psi$. If $\beta$ is a reflection, then 
$\langle \beta^\Delta\rangle{\,=\,}\Psi$ is a compact group 
of homologies by \cite{cp} 61.2, and $\dim\Psi{\,=\,}3$. \\
(c) Let $\beta$ be a \emph{planar} involution in $\Psi$,
and let $\Upsilon$ be a semi-simple subgroup of the 
centralizer $\Cs\Delta\beta$. Then 
$\dim\Upsilon{\,\le\,}13$: Write 
$\Upsilon|_{\cF_\beta}{\,=\,}\Upsilon/\Kappa$ 
and $\Lambda{\,=\,}\Kappa^1$. The kernel $\Lambda$ is 
compact and semi-simple, hence a Lie group, and then
$\Lambda$ is contained in $\Spin3\RR$ by stiffness. 
In \cite{65} 6.1 it has been shown that 
$\Upsilon{:}\Kappa{\,\le\,}10$. \\
(d) If $\dim\Psi{\,=\,}16$, then 
$\Psi{\,\cong\,}{\rm(P)}\SL3\CC$. The centralizer of 
the involution  $\beta{\,=\,}{\rm diag}(1,-1,-1)$ contains a group
$\SL2\CC{\cdot}\Gamma$ of dimension $6{+}9$ and is too 
large. \\
(e) Each group $\Psi$ of type ${\rm A}_3$ contains a planar involution $\beta$ centralizing 
a $6$-dimensional semi-simple subgroup of $\Psi$. This is easy in the cases $\SU4\CC$,\,
$\SL2\HH$,\, $\SL4\RR$ and its universal covering group; 
it is less obvious for the groups related to $\SU4(\CC,r)$, $r{\,=\,}1,2$, which are not necessarily Lie groups.
In both cases, a suitable involution $\beta$ is determined 
by ${\rm diag}(-1,-1,1,1)$; then $\Cs\Psi\beta$ 
contains a group  locally isomorphic to 
$\SU2\CC{\times}\SL2\RR$ or to  
$(\SU2\CC)^2$, respectively, and  
$\Cs\Delta\beta$ is too large. \\
(f) Each involution of the compact group $\Gtwo$ 
is centralized by $\SO4\RR$, see \cite{cp} 11.31. 
The non-compact groups of type $\Gtwo$ contain 
either $\SO4\RR$ or $\Spin4\RR$. In any case, the 
centralizer of a non-central involution is too large. \\
(g)  If $\Psi$ is locally isomorphic to a group 
$\U2(\HH,r)$ of dimension $10$, then 
$\beta{\,=\,}{\rm diag(1,-1)}$ describes a \emph{planar} 
involution centralized by a semi-simple group of dimension
$6{+}15$, which is too large. \\
(h) If $\dim\Psi{\,=\,}10$, then $\Psi$ maps onto 
$\Sp4\RR$ or $\Opr5(\RR,2)$. In the latter case, the 
diagonal element ${\rm diag(- \1,1)}$ corresponds to a \emph{planar} 
involution $\beta$ such that 
$\Cs\Psi\beta{\,\circeq\,}\Opr4(\RR,1){\,\cong\,}\SO3\CC$, again a contradiction. \\
(i) Finally, let 
$\Psi/(\Zeta\smcap\Psi){\,\cong\,}\Sp4\RR$.  
Conceivably, $\Psi$ is not a Lie group. Again $\Psi$ 
cannot contain a planar involution. Hence the unique 
involution in $\Psi$ is a reflection 
$\sigma{\,\in\,}\Psi_{[u,av]}{\,:=\,}\Eta$. 
As $\dim\Eta{\,\le\,}8$, there is a point $p$ such that 
$\cG{\,=\,}\langle p^{\Psi}\rangle$ is a connected 
subplane,\hskip3pt and $\cG{\hskip2pt\ledot\,}\cP$\hskip-3pt, or else $\dim\Psi|_\cG{\,\le\,}6$.  
We have $\Gamma_{\hskip-2pt p}|_\cG{\,=\,}\1$, and 
$\Gamma_{\hskip-2pt p}$ is a compact Lie group by step (a). 
Stiffness 2.5(\^b) shows $\dim\Gamma_{\hskip-2pt p}{\,\le\,}3$ and 
$\dim p^\Gamma{\,\ge\,}12$. Consequently 
$\langle p^\Gamma\rangle{\,=\,}\cP$ and 
$\Psi_{\hskip-3pt p}{\,=\,}\1$. It follows that 
$\cG{\,=\,}\cP$ and $\Gamma_{\hskip-2pt p}{\,=\,}\1$. 
Split $\Gamma$ into an almost direct product $\Alpha\Beta$ with 
$\dim\Alpha{\,=\,}6$.\hskip3pt Then 
$\dim p^\Beta{>\,}8$,\; 
$\langle p^\Beta\rangle{=\,}\cP$, and 
$(\Psi\Alpha)_p{\,=\,}\1$.
Hence $p^{\Psi\Alpha}$ is open in $P$ by 
\cite{cp} 53.1),  and $p^\Delta$ is also open. This implies that $\Delta$ is a Lie group, see 2.6. \\
(j) The cases $\dim\Psi{\,=\,}8$ can be treated similarly. 
Again there is a point $p$ such that 
$\cG{\,=\,}\langle p^\Psi\rangle{\,\ledot\,}\cP$, 
$\Gamma_{\hskip-2pt p}|_\cG{\,=\,}\1$,\; 
$\dim p^\Gamma{\ge\,}14$,\;  
$\langle p^\Gamma\rangle{\,=\,}\cP$, and 
$\Psi_p{\,=\,}\1$. 
Let $\Beta$ be the other $8$-dimensional factor of $\Delta$ and denote the product of the remaining factors by 
$\Upsilon$. Then 
$\langle p^{\Beta\Psi}\rangle{\,=\,}\cP$, 
or else an action of $\Beta\Psi$ on $\cG$ 
would have a kernel of dimension $0$.  Now 
$\Upsilon_{\hskip-3pt p}{\,=\,}\1$,\; 
$\langle p^\Upsilon\rangle{\,=\,}\cP$,\; 
$(\Beta\Psi)_p{\,=\,}\1$,\; $p^\Delta$ is open in 
$P$, and $\Delta$ is a Lie group by   2.6. \Qed
\par\bigskip
{\Bf 8. Fixed triangle}
\par\medskip
{\bold 8.1.}
{\it If $\cF_\Delta{\,=\,}\langle o,u,v\rangle$ is a triangle, 
if $\Delta$ is semi-simple, and if $\dim\Delta{\,\ge\,}21$, 
then $\Delta$ is a Lie group\/}. 
\par\smallskip
{\tt Proof.} For groups of dimension ${\ge\,}22$, the claim is true by \cite{66} 9.2 Corollary. Hence we may suppose that 
$\dim\Delta{\,=\,}21$. Then $\Delta$ is almost simple or 
one of the following holds: ({\it i\/}) $\Delta$ has a $15$-dimen\-sional factor, ({\it ii\/}) $\Delta$ is a product of $3$ factors of dimensions $10,8,3$, ({\it iii\/}) each factor of $\Delta$ has dimension $3$ or $6$. Let again 
$\Zeta{\,=\,}\Cs{}\Delta$ be the center of $\Delta$ and suppose that $\Zeta$ is not a Lie group. Note 
that $\Zeta$ acts freely on the complement $D$ of
$ou{\,\smcup\,}ov{\,\smcup\,}uv$; in fact, 
$p^\zeta{=\,}p$,\; $\1{\,\ne\,}\zeta{\hskip1pt\in\hskip1pt}\Zeta$ 
implies $\cD{\,=\,}\langle p^\Delta,o,u,v\rangle{\,<\,}\cP$. 
Put $\Delta|_\cD{\,=\,}\Delta/\Kappa$. Then 
$\Delta{:}\Kappa{\,\le\,}11$\; (see \cite{cp} 83.26),\; 
$\dim\Kappa{\,\ge\,}10$,\; 
$\cD$ is flat by stiffness, $\Delta|_\cD$ is trivial, and $\Kappa$ is too large. \\
(a) If $\Delta$ is almost simple, then $\Delta$ is a Lie 
group, or $\overline\Delta{\,=\,}\Delta/\Zeta$ is isomorphic to  $\Opr7(\RR,2)$ or to $\PSp6\RR$. \\
(b) In the first case, $\overline\Delta$ has a subgroup 
$\SO5\RR$; it can be lifted to a group 
{\cyss Yu}${\,\cong\,}\U2\HH{\,\le\,}\Delta$, see \cite{cp} 94.27 or 2.9, and note that  by 2.4 there is no 
subgroup $\SO5\RR$ in $\Delta$. 
The central involution $\iota$ of {\cyss Yu} is 
not planar (or else 
{\cyss Yu}$|_{\cF_\iota}{\,\cong\,}\SO5\RR$ contrary to 2.4).
Therefore $\iota$ is a reflection; it   is in the kernel of the action of $\Delta$ on one side of the fixed triangle and hence $\iota{\,\in\,}\Zeta$; but then 
$\iota{\,=\,}\overline\iota$ would be trivial. \\ 
(c) In the second case, the subgroup 
$\SU3\CC$ of $\overline\Delta$ can be lifted to an
isomorphic copy $\Omega{\,<\,}\Delta$, and $\Omega$ 
is contained in a maximal compact subgroup 
$\Xi{\,<\,}\Delta$; moreover $\dim\Xi{\,=\,}7$ and  
$\Xi{\,=\,}\Theta\Omega$ is connected, $\Theta$ is 
commutative and $\Theta{\,\le\,}\Cs{}\Omega$\, (\cite{cp} 
93.10,11).  Let $x$ be an arbitrary point in $D$\, (the 
complement of the fixed lines). Then $\Zeta$ acts freely 
on the subplane $\cD{\,=\,}\cF_{\hskip-2pt\Delta_x}$, and $\Delta_x$ is a Lie group because 
$\Delta_x\smcap\Zeta{\,=\,}\1$. 
If $\Delta$ has an $8$-dimensional orbit on two sides of 
the fixed triangle, then $\Delta$ induces a Lie group 
on these orbits (see 2.6), and $\Delta$ itself would be a 
Lie group. Consequently  $\dim x^\Delta{\,\le\,}14$ and $\dim\Delta_x{\ge\,}7$.
As $\Zeta$ is not a Lie group, 2.12 implies 
$\dim\cD{\,\in\,}\{0,8\}$. In the second case,  
Stiffness 2.5(\^b) shows that  $\Delta_x{\,\cong\,}\Spin3\RR$.
Hence $\dim\cD{\,=\,}0$ for each $x{\,\in\,}D$. 
There are $3$ pairwise commuting conjugate involutions 
in $\Omega$; these are \emph{planar} because axes of 
reflections are fixed by $\Delta$. An involution 
$\iota{\,\in\,}\Omega$   is centralized in $\Omega$  by a group $\Phi{\,\cong\,}\SU2\CC$. In the case 
$\Phi|_{\cF_\iota}{\,\ne\,}\1$ there would exist 
\emph{planar} involutions on $\cF_\iota$, and $\Zeta$ 
would be a Lie group. Hence $\cF_\iota{\,=\,}\cF_\Phi$.
All involutions  $\kappa{\,\in\,}\Omega$ such that 
$\iota\kappa{\,=\,}\kappa\iota$ generate a subgroup 
$\Psi{\,\cong\,}\SU2\CC$. As there are no $4$-dimensional 
subplanes, $\kappa$ and its conjugates induce reflections on  $\cF_\iota$. Therfore $\Psi|_{\cF_\iota}$ is a group of homologies, say with center $o$. \\
On the other hand, 
there is a product $\Gamma\Upsilon{\,<\,}\Delta$ corresponding to $\Sp2\RR{\times}\Sp4\RR$ in a double covering of $\overline\Delta$\, (use \cite{cp} 94.27). 
Note that $\Sp2\RR{\,=\,}\SL2\RR$ and that $\Sp4\RR$ has a subgroup $\SU2\CC$. Suppose that the latter group is mapped onto a group $\SO3\RR$ in $\Upsilon$. The arguments above show that $\Upsilon$ contains a circle group of homologies, but this is impossible, since $\Upsilon$
is almost simple. Therefore $\Upsilon$ contains a central 
involution.  Up to conjugacy, 
$\iota{\,\in\,}\Cs{}{\hskip-3pt\Upsilon}$, and 
$\dim(\Theta\Upsilon)|_{\cF_\iota}{\,\ge\,}10$.
According to \cite{69} 8.0, it follows that  
$(\Cs\Delta\iota)^1|_{\cF_\iota}$ is a Lie group. 
Recall that $\Zeta$ acts freely on $D$. Hence $\Zeta$ 
is a Lie group, and so is $\Delta$  \\
({\it i\/}) A $15$-dimensional factor $\Upsilon$ is a Lie group or $\Upsilon/\Zeta$ is locally isomorphic to 
$\SU4(\CC,r), r=1,2$. The group $\Delta$ is an almost 
direct product $\Upsilon\Gamma$. By stiffness, 
$p^{\hskip-1.5pt\Upsilon}{\ne\,}p$ for each $p{\,\in\,}D$, and then $\langle p^{\hskip-2pt\Upsilon\!},o,u,v\rangle{\,=\,}\cP$. Consequently, $\Gamma$ acts freely on $D$ and 
$\langle p^\Gamma\!,o,u,v\rangle{\,\ledot\,}\cP$ by 2.12.  \\
(d) Suppose in case ({\it i\/}) that $r{\,=\,}1$.  Then 
$\Upsilon/\Zeta$ has a subgroup $\SU3\CC$. 
As in step (c)  there is a group $\Omega{\,\cong\,}\SU3\CC$  
in $\Upsilon$, and all involutions in $\Omega$ are 
\emph{planar}. Choose some involution 
$\iota{\,\in\,}\Omega$. Using the arguments and notation 
of step~(c), we obtain the following:
$\SU2\CC{\,\cong\,}\Phi{\,\le\,}\Cs\Omega\iota$,
$\cF_{\hskip-2pt\iota}{=\,}\cF_\Phi$,
$\SU2\CC{\,\cong\,}\Psi{\,\le\,}\Cs\Omega\iota$,
and  $\Psi|_{\cF_\iota}$ is a group of homologies.
The center $\Zeta$ acts freely on $D\smcap\cF_\iota$.
Hence the kernel of $\Gamma|_{\cF_\Phi}$ is a compact 
Lie group; by stiffness $\dim\Gamma|_{\cF_\iota}{\,=\,}6$,
and $\Theta\Psi\Gamma$ induces a $10$-dimensional 
group on $\cF_\iota$. According to \cite{69} 8.0
the induced group is a Lie group. Therefore $\Zeta$ and 
$\Delta$ are Lie groups. \\
(e) If $r{\,=\,}2$, then 
$\Upsilon/\Zeta{\,\cong\,}\Opr6(\RR,2)$ has a subgroup 
$\SO4\RR$. By 2.9 there is a subgroup 
$\Chi{\,<\,}\Upsilon$ which is isomorphic to
$\SO4\RR$ or to $\Spin4\RR$. If 
$\Chi$ contains a reflection, say
$\iota{\,\in\,}\Delta_{[o,uv]}$, then $\iota{\,\in\,}\Zeta$\, 
(or else $\langle \iota^{\hskip-2pt\Upsilon}\rangle{\,=\,}\Upsilon{\,\le\,}\Delta_{[o,uv]}$, which is impossible). 
Recall that $\Gamma$ acts freely on $D$. Suppose first 
that $\Chi{\,\cong\,}\SO4\RR$. Then all involutions in $\Chi$ are \emph{planar}. Let $\beta$ denote the central involution in $\Chi$. There is a subgroup $\Phi{\,\cong\,}\SO3\RR$ 
in $\Chi$. The kernel of $\Chi|_{\cF_\beta}$ is a compact Lie group;  by stiffness it is contained in $\Spin3\RR$. Therefore 
$\Phi$ acts faithfully on $\cF_\beta$, and 
$\cC{\,=\,}\cF_{\alpha,\beta}{\,\ldot}\cF_\beta$ for each involution $\alpha{\,\in\,}\Phi$, but  
$\cC^\Gamma{=\,}\cC$, and orbits of $\Gamma$ on
$D\smcap\cC$ are $6$-dimensional, a contradiction.
Consequently 
$\Chi{\,\cong\,}\Spin4\RR{\,=\,}(\Spin3\RR)^2$.
Only one of the $3$ involutions in $\Chi$ can be contained in the center $\Zeta$. Let $\beta$ be a \emph{planar} 
involutin in $\Chi$ such that 
$\Chi/\langle\beta\rangle{\,\cong\,}\SO3\RR{\times}\Spin3\RR$. If $\iota$ is any involution in the facor 
 $\Phi{\,\cong\,}\SO3\RR$, then 
 $\iota|_{\cF_\beta}{\,=\,}\1$ and, as $\Phi$ is simple, 
 even $\Phi|_{\cF_\beta}{\,=\,}\1$, but this contradicts stiffness  2.5(\^b). \\
({\it ii\/}) $\Delta$ has $3$ almost simple factors 
$\Upsilon,\Psi,\Gamma$ of dimensions $10,8,3$, respectively. Recall that a semi-simple group of an $8$-dimensional plane with a triangle of fixed elements has dimension at most $9$\, (\cite{63} 6.1 = ($\dagger$)\,). By stiffness it follows that $\langle p^\Delta,o,u,v\rangle{\,=\,}\cP$. If $p^\Upsilon{\,=\,}p$, then
$\langle p^{\Psi\Gamma}\rangle{\,=\,}\cP$, 
$\Upsilon_{\hskip-2pt p}{\,=\,}\1$, and 
$p^{\hskip-2pt\Upsilon}{\ne\,}p$. By ($\dagger$)
we have $\langle p^\Upsilon\rangle{\,=\,}\cP$.
If  $\cD{\,=\,}\langle p^{\Psi\Gamma}\rangle{\,\ldot}\cP$, then 
($\dagger$) implies $\Gamma|_\cD{\,=\,}\1$, and 
$\Gamma$ is a compact Lie group. Stiffness shows  
$\Gamma{\,\cong\,}\Spin3\RR$. The involution 
$\iota{\,\in\,}\Gamma$ is planar, and 
$\cF_{\hskip-2pt\iota}^\Delta{\,=\,}\cF_\iota$, a contradiction. Therefore $\langle p^{\Psi\Gamma}\rangle{\,=\,}\cP$, and hence 
$\Upsilon_{\hskip-3pt p}{\,=\,}\Psi_{\hskip-2pt p}{\,=\,}\1$.
 The different possibilities will be treated separately. \\
(f) If $\Upsilon$ is the compact group $\U2\HH$, then
each non-central involution $\iota$ in $\Upsilon$ is planar, since $\langle \iota^{\!\Upsilon}\rangle{\,=\,}\Upsilon$ 
is not axial. But then $\iota{\,\in\,}\Upsilon_{\!p}$ for each fixed point $p$ of $\iota$. This contradicts the statements at the end of ({\it ii\/}). Analogous arguments show that 
$\Psi$ is not compact. \\
(g) If $\Upsilon$ is locally isomorphic to 
$\U2(\HH,1)$, then $\Upsilon$ is simply connected, 
because $\Opr5(\RR,1)$ has a subgroup 
$\Phi{\,\cong\,}\SO3\RR$ and $\cF_\Phi$ is flat. 
Therefore $\Upsilon$ has a maximal compact subgroup 
$\Chi{\,\cong\,}\Spin4\RR$. Again a non-central involution 
$\beta$ is planar, which contradicts ({\it ii\/}). \\
(h) Only the possibility 
$\Upsilon/\Zeta{\,\cong\,}{\rm(P)}\Sp4\RR$ remains. 
Again $\Delta$ has no subgroup $\SO3\RR$. 
Either $\Psi/\Zeta$ is locally isomorphic to 
$\SU3(\CC,1)$ or $\Psi$ is a double covering $\Pi$ of 
$\SL3\RR$. 
A subgroup $\SU2\CC$ of the unitary group 
can be lifted to an isomorphic copy 
$\Phi{\,<\,}\Psi$; the involution $\iota{\,\in\,}\Phi$ is not 
central and hence it is planar, which contradicts ({\it ii\/}). 
Thus $\Psi{\,\cong\,}\Pi$. Both factors $\Upsilon$ and 
$\Psi$ have a subgroup $\Spin3\RR$. The central involutions  of these subgroups are reflections (because 
$\Upsilon_{\hskip-3pt p}{\,=\,}\Psi_{\hskip-2pt p}{\,=\,}\1$), and these are distinct, or else there would exist a subgroup 
$\Omega{\,\cong\,}\SO3\RR$. 
A maximal compact subgroup of $\Upsilon\Psi$ is at most $7$-dimensional.  It contains a subgroup 
$\Chi{\,\cong\,}\Spin4\RR$ without any \emph{planar} involution. \\
(i) {\it $\cG{\,=\,}\langle p^\Gamma\!,o,u,v\rangle{\,\ldot}\cP$ for any $p{\,\in\,}D$\/}: if $p^\Gamma{=\,}p$, then 
$p^\Delta{\,\subseteq\,}\cF_\Gamma{\,<\,}\cP$ and 
$\dim \Delta_p{\,\ge\,}13$. Stiffness would imply 
$\Delta_p^1{\,\cong\,}\Gtwo$, which is not true.
Hence $\Gamma$ acts almost effectively on $\cG$ and 
$\cG{\,\ledot\,}\cP$.\hskip2pt Equality is impossible, since 
$\dim(\Upsilon\Psi)_p{\,>\,}1$. Denote the connected component of $(\Upsilon\Psi)_p$ by $\Lambda$.
By 2.6 we may assume that 
$\dim\Lambda{\,\ge\,}3$. Note that 
$\Lambda|_\cG{\,=\,}\1$. Hence $\Lambda$ is compact.
It follows that $\dim(\Lambda\smcap\Chi){\,\ge\,}2$.  
Obviously, each involution in \
$\Lambda\smcap\Chi$ is planar. This contradicts the last statement in step (h). \\
({\it iii\/}) Each $6$-dimensional factor of $\Delta$ is locally 
and then even globally isomorpic to $\SL2\CC$.\hskip2pt  There is at least one factor $\Gamma$ such that  
$\dim\Gamma{=\,}3$.  We may assume that $\Gamma$ is not a Lie group. Then a maximal compact subgroup  $\Theta$ of 
$\Gamma$  is a $1$-dimensional  connected subgroup. Let  
$\Xi $ denote the product of the other factors of $\Delta$ and apply  \cite{66} 9.2 (B) to the group $\Eta{\,=\,}\Xi\Theta$. 
As $\Theta$ is normal in $\Eta$, we obtain 
$\dim\Eta{\,<\,}18$ and  $\dim\Delta{\,\le\,}19$. \Qed
\par\bigskip
 
{\Bf  Summary}
\par\medskip
{\it If $\Theta$ denotes a compact connected $1$-dimensional subgroup  and
if $\dim\Delta{\,\ge\,}k$, the number in the table, then $\Delta$ is a Lie 
group\/}. Bold face entries improve corresponding results in \cite{66}. 

\par\medskip
 \begin{center}
\begin{tabular}{|c||c|c|c|c||l|} 
\hline
$\cF_\Delta$ & $\Delta$ $s$-$s$  & $\Theta\triangleleft \Delta$ & 
$\RR^t\triangleleft\Delta$ & $\Delta$ arbitr. & References   \\  \hline 
$\emptyset$ & {\bf 20} &  && {\bf 22} &  3.1, 3.2  \\
$\{W\}$ & {\bf 11} &  && {\bf 22} & 4.2, 4.0 \\
flag & {\bf 11} &  && 22 & 5.2, \cite{66}\,5.0 \\  \hline
$\{o,W\}$ &\hskip12pt {\bf 25}\,$^{1)}$& && 27 & 
6.4), \cite{psz} \\  
$\langle u,v\rangle$ & {\bf 11} & && 18 & 5.2, 
\cite{66}\,6.0 \\
$\langle u,v,w\rangle$ & {\bf 11} & && 15 & 5.2, 
\cite{66}\,7.0 \\  \hline
$\langle u,v,ov\rangle$ & {\bf 25} & $20\,^{2)}$ && 27 & 7.1, \cite{psz}  \\  
$\langle o,u,v\rangle$ & {\bf 21} & $18\,^{\ \ }$  & 22 & 22 & 8.1, \cite{66}\,9.2 \\
arbitrary & {\bf 25} &  && 27&  \cite{psz} \\  \hline
\end{tabular}
\end{center}
\par 
\begin{tabular}{l l l l}
\quad $^{1)}$ & $\cF_\zeta{\,<\,}\cP{\,\Rightarrow\,}\dim\Delta{\,\le\,}19$\; (6.1)  \cr
\quad $^{2)}$ & or $\Cs{}\Delta|_{av}{\,=\,}\1$\; (8.3 Remark) \cr
\end{tabular}
\par

\par
\par\bigskip
\cite{66} {\bf 9.2.} {\it Assume that $\Delta$ fixes a triangle $a,u,v$\/}. \\
(A) {\it If $\Delta$ is semi-simple and if $\dim\Delta{\,\ge\,}22$, then $\Delta$ is a Lie group\/}. \\
(B) {\it If $\Delta$ has a $1$-dimensional compact normal subgroup $\Theta$  and if 
 $\dim\Delta{\,\ge\,}18$, then $\Delta$ is a Lie group and $\dim\Delta{\,\le\,}23$\/}. \\
(C) {\it If $\Delta$ has a minimal normal vector subgroup $\Theta{\,\cong\,}\RR^t$  and if 
$\dim\Delta{\,\ge\,}22$, then $\Delta$ is a Lie group\/}. 
\par\smallskip
({\tt Proof} of \cite{66}\,9.2(B). All orbits on two sides of the fixed triangle have dimension~${<\,}8$, and $\dim p^\Eta{\,\le\,}14$ for each 
point $p{\,\notin\,}au{\smcup}av{\smcup}uv$. 
Put $\Lambda{\,=\,}(\Eta_p)^1$. By the approximation theorem \cite{cp} 93.8, there is a compact $0$-dimensional central subgroup $\Nu{\,\triangleleft\,}\Eta$ such that 
$\Eta/\Nu$ is a Lie group. Note that $\Theta$ is in the center of $\Eta$. If $p^\Theta{\,=\,}p$, then 
$\cF_\Theta{\,=\,}\cF_{\hskip-1pt\Theta}^{\:\Eta}{\,<\,}\cP$ and 
$\Eta{\hskip1pt:\hskip1pt}\Eta_p{=}\dim p^\Eta{\le}\dim\cF_\Theta{=}d$ with  $d{\,\in\,}\{0,2,4,8\}$. Either 
$d{\le}2$,\,  $\Theta{\,\triangleleft\,}\Lambda{\,=\,}(\Eta_p)^1{\,\not\cong\,}\Gtwo$,\, 
$\dim\Lambda{\,\le\,}11$,  and $\dim\Eta{\,\le\,}13$, or 
$d{\,\ge\,}4$,  $\dim\Lambda{\,\le\,}8$, 
and $\dim\Eta{\,\le\,}16$. If $p^\Theta{\,\ne\,}p$, then $\Theta$ acts non-trivially on the subplane 
$\cE{\,=\,}\langle p^{\Theta\Nu},a,u,v \rangle$, $\cE$ is not flat by \cite{cp} 32.17, and 
$\dim\cE{\,=\,}4$ or $\cE{\,\ledot\,}\cP$. In the first case, $\Nu|_\cE{\,=\,}\Nu/\Kappa$ is a Lie group by \cite{cp} 71.2, and $\Kappa{\,\ne\,}\1$. Hence 
$\cF_\Kappa{\,=\,}\cF_\Kappa^{\hskip2pt\Eta}{\,<\,}\cP$,\, $p^\Eta{\,\subseteq\,}\cF_\Kappa$,\, $\Eta{\hskip1pt:\hskip1pt}\Eta_p{\,\le\,}8$, 
$\Eta_p|_\cE{\,=\,}\1$,\, $\dim\Eta_p{\,\le\,}8$, and 
$\dim\Eta{\,\le\,}16$. In the second case, we have 
$\Nu_p|_{\langle p^\Eta\rangle}{\,=\,}\1$. Therefore 
$\Eta{\hskip1pt:\hskip1pt}\Eta_p{\,=\,}\dim p^\Eta{\,\le\,}8$ or 
$\Eta_p\smcap\Nu{\,=\,}\1$ and $\Eta_p$ is a Lie group. Stiffness implies $\dim\Eta_p{\,\le\,}7$ or 
$\dim\Eta_p{\,\le\,}3$, respectively, and 
$\dim\Eta{\,\le\,}17$.) \\

\qquad Helmut R. Salzmann \par
\qquad Mathematisches Institut \par
\qquad Auf der Morgenstelle 10 \par
\qquad D-72076 T\"ubingen \par
\vspace{5pt}
\qquad helmut.salzmann@uni-tuebingen.de  
\end{document}